\documentclass{amsart}
\usepackage{amssymb,eepic,epic}
\newtheorem{theorem}{Theorem}[section]
\newtheorem{corollary}[theorem]{Corollary}
\newtheorem{proposition}[theorem]{Proposition}

\newtheorem{lemma}[theorem]{Lemma}
\theoremstyle{definition}
\newtheorem{definition}[theorem]{Definition}
\theoremstyle{remark}
\newtheorem{remark}[theorem]{Remark}
\newtheorem{example}[theorem]{Example}
\newcommand\A{\mathcal{A}}
\newcommand\M{\mathcal{M}}
\newcommand\G{\mathcal{G}}

\renewcommand{\O}{\mathcal{O}}

\newcommand{\N}{\mathbb{N}}
\newcommand{\R}{\mathbb{R}}
\newcommand{\C}{\mathbb{C}}
\newcommand{\Z}{\mathbb{Z}}
\newcommand{\Q}{\mathbb{Q}}

\newcommand\lie[1]{\mathfrak{#1}}
\renewcommand{\k}{\lie{k}}
\newcommand{\h}{\lie{h}}
\newcommand{\g}{\lie{g}}

\newcommand{\z}{\lie{z}}
\renewcommand{\t}{\lie{t}}
\newcommand{\Alc}{\lie{A}}

\newcommand{\RR}{ \operatorname {RR} }
 
\newcommand{\Ad}{ \operatorname {Ad} } 
\newcommand{\Hol}{ \operatorname {Hol} } 
\newcommand{\Rep}{\operatorname{Rep}}

\newcommand{\Ind}{ \operatorname {Ind}}
\renewcommand{\ker}{ \operatorname {ker}}
\newcommand{\im}{ \operatorname {im}}

\newcommand{\Spin}{ \operatorname{Spin}}

\newcommand\dirac{/\kern-1.2ex\partial} % Dirac operator
\newcommand\qu{/\kern-.7ex/} % Categorical quotients
\newcommand{\fus}{\circledast}
\newcommand{\lev}{{\lambda}} % Level
\newcommand{\levi}{{m}} % Integral Level

\newcommand{\hra}{\hookrightarrow}
\newcommand{\ra}{\rightarrow}
\newcommand{\bib}{\bibitem}

\renewcommand{\d}{{\mbox{d}}}
\newcommand{\ol}{\overline}
\newcommand\Gpar{\G_\partial}
\newcommand\Phinv{\Phi^{-1}}

\newcommand\Sig{\Sigma}
\newcommand\sig{\sigma}
\newcommand\eps{\epsilon}

\newcommand\om{\omega}

\newcommand{\del}{\delta}
\newcommand{\f}{\frac}

\newcommand{\p}{\partial}
\renewcommand{\l}{\langle}
\renewcommand{\r}{\rangle}

\newcommand{\ti}{\tilde}

\begin{document}

%\sloppy

\title{A symplectic proof of Verlinde factorization}

\author{E. Meinrenken}
\address{Massachusetts Institute of Technology, Department
of Mathematics, Cambridge, Massachusetts 02139}
\email{mein@math.mit.edu}

\author{C. Woodward}
\address{Harvard University, Department of Mathematics, 1 Oxford Street,
Cambridge, Massachusetts 02138}
\email{ woodward@math.harvard.edu}

\begin{abstract}
We prove a multiplicity formula for Riemann-Roch numbers of reductions
of Hamiltonian actions of loop groups.  This 
includes as a special case the 
factorization formula for the quantum dimension of the moduli space of
flat connections over a Riemann surface.\\  
\end{abstract}

\maketitle

\tableofcontents

\newpage

\section{Introduction}

The quantization of the moduli space of flat connections over a
Riemann surface has been the subject of intensive study from a number
of different points of view.  Much of the recent work in mathematics
has focused on proving formulas for the dimension of the quantization
discovered by the physicist E. Verlinde \cite{V} in the context of
conformal field theory.  Essentially there are two ingredients in
Verlinde's approach: the ``factorization theorem'' (included by Segal
as one of the axioms of conformal field theory \cite{Seg}) describes
the quantization of the moduli space associated to a Riemann surface
obtained by gluing a second (possibly disconnected) surface along two
boundary circles; the ``fusion rules'' describe the quantization of
the moduli space of a three-holed sphere with boundary components
marked by irreducible representations of the loop group.  Together
these give a formula for the dimension of the quantization in terms of
a pants decomposition of the surface.  Both parts of Verlinde's
approach were carried out rigorously by Tsuchiya-Ueno-Yamada
\cite{TUY}.  Since then there has been a vast amount of work on
improving and understanding the formulas: see in particular the works
by Beauville-Laszlo \cite{BL}, Bertram-Szenes \cite{BS},
Daskalopoulos-Wentworth \cite{DW}, Faltings \cite{F},
Kumar-Narasimhan-Ramanathan \cite{KNR}, Teleman \cite{Tel}, and
Thaddeus \cite{Tha}, among others.  There is an alternative approach,
initiated by Witten \cite{W1} and carried out by Szenes \cite{Sze} and
Jeffrey-Kirwan \cite{JK}, in which the Verlinde formulas are derived
from the cohomology ring of the moduli space.  Very recently, yet
another approach was outlined by S. Chang \cite{C3}, who obtained a
character formula for quantizations of Hamiltonian loop group actions
and showed that this implies Verlinde's formula.

In this paper, we derive the factorization theorem in the context of
symplectic geometry.  To explain the idea (which has appeared in the
literature in various forms) let $\Sig$ be any compact oriented
Riemann surface with boundary, and $\M(\Sig)$ the space of flat
$G$-connections on $\Sig$ mod gauge transformations that are trivial
on the boundary.  The moduli space $\M(\Sig)$ carries a Hamiltonian
action of $b$ copies of the loop group $LG$ where $b$ is the number of
boundary components of $\Sig$.  If $\Sig$ is formed from a second
surface $\hat{\Sig}$ by gluing together two boundary components, then
$\M(\Sig)$ is related to $\M(\hat{\Sig})$ by a symplectic reduction by
$LG$.  Assuming heuristically that some quantization procedure
constructs a projective $LG^b$-representation $Q(\M(\hat{\Sig}))$ of
$\M(\hat{\Sig})$, the ``quantization commutes with reduction''
principle \cite{GS2} implies that $Q(\M(\Sig))$ should be the
$LG$-invariant part of $Q(\M(\hat{\Sig}))$.  A corollary of this
principle is a ``factorization formula'' for the quantizations of
finite-dimensional symplectic quotients of $\M(\Sig)$ in terms of
those of $\M(\hat{\Sig})$.

Our main result Theorem \ref{MainResult} is a formula of this type in
the general setting of Hamiltonian loop group actions with proper
moment maps.  The proof is entirely finite-dimensional.  Because of
the properness assumption, the symplectic quotients are
finite-dimensional, and their quantizations can be defined as indices
of suitable Dirac operators (Riemann-Roch numbers), using
desingularizations if necessary.  Using finite-dimensional
``symplectic cross-sections'' for Hamiltonian loop group actions we
reduce the proof to the finite-dimensional version of ``quantization
commutes with reduction'', which has been proved in general in
Meinrenken, Meinrenken-Sjamaar \cite{M2,MS}.  An important ingredient
is a gluing formula for the behavior of Riemann-Roch numbers under
``symplectic surgery'' \cite{M2,MW}.

Many of the ideas used in our proof are already present in the
literature.  The construction of the moduli space of a Riemann surface
without boundary as a symplectic quotient of the space of all
connections by the gauge group action goes back to Atiyah and Bott
\cite{AB}.  For Riemann surfaces with non-empty boundary the
corresponding moduli space was constructed by Donaldson \cite{D} as a
special case of moduli spaces of framed connections.  The relationship
between the Verlinde formula and ``quantization commutes with
reduction'' was outlined by Segal \cite{Seg2} and is present both in
the algebraic geometry approach and also in the physics literature
(e.g. \cite{W1,EMSS}).  The symplectic cross-sections used here are
generalizations of the ``extended moduli space'' of Chang \cite{C1},
Huebschmann \cite{H} and Jeffrey \cite{J}. Some of our results on
Hamiltonian $LG$-actions overlap with those of Chang \cite{C2}.  The
densely-defined torus actions in this paper are related to the work of
Goldman \cite{G} and Jeffrey-Weitsman \cite{JW} on the ``twist flows''
on the moduli space.\\

\noindent {\bf Acknowledgments} We thank S. Chang, L. Jeffrey,
S. Martin, A. Szenes, and C. Teleman for helpful discussions.

\section{Statement of results}
\label{SubsectionGluingEqualsReduction}

In this section we describe our approach to the factorization
property.
In order to state the main result as quickly as possible we omit the
details of Sobolev topologies and existence of infinite-dimensional
quotients which are postponed until Appendix \ref{YangMills}.

\subsection{Construction of the moduli space}
\label{ConstructionModuliSpace}
We begin by reviewing Yang-Mills theory for Riemann surfaces with
boundary along the lines of Atiyah-Bott \cite{AB} and Donaldson
\cite{D}.  Let $\Sigma$ be an oriented compact Riemann surface with
$b$ boundary components. If $\Sigma$ is connected and has genus $g$,
we will write $\Sigma =\Sigma^b_g$.  We denote by
$\iota:\,\partial\Sigma\hra \Sigma$ the inclusion of the boundary.

Let $G$ be a connected and simply-connected compact Lie group. We fix
an $\text{Ad}$-invariant inner product on the Lie algebra $\g$,
normalized by the requirement that on each simple summand of $\g$, the
norm squared of the coroot corresponding to the highest root is equal
to $2$.  Let $\mathcal{P}$ be a principal $G$-bundle over $\Sigma$,
and $\G(\Sigma)$ the group of gauge transformations, i.e. the space of
sections of the associated bundle $\mathcal{P}\times_{Ad(G)} G$.
Since $\pi_0(G) = \pi_1(G) = \{ 1 \}$, the bundle $\mathcal{P}$ is
necessarily trivial and we can identify $\G(\Sigma)$ with the space of
maps $\Sigma\ra G$ and its Lie algebra with $\Omega^0(\Sigma,\g)$.
Let $\A(\Sigma)\cong \Omega^1(\Sigma,\g)$ be the space of principal
connections on $\mathcal{P}$.  It has a symplectic form given by
$$ \omega_A(a_1,a_2)=\int_\Sigma a_1 \stackrel{ _\cdot}{\wedge} a_2 $$
making
$\A(\Sig)$ into an infinite-dimensional symplectic manifold.  For all
$A\in\A(\Sigma)$, we denote by
$$\d_A=\d+[A,\cdot]: \Omega^i(\Sigma,\g)\ra
\Omega^{i+1}(\Sigma,\g)$$ 
the associated covariant derivative and by $F_A=\d\,A+\f{1}{2}[A,A]\in
\Omega^{2}(\Sigma,\g)$ its curvature.  The natural gauge group action
$$ g\cdot A = \text{Ad}_g(A)- \d g\, g^{-1}$$ preserves the symplectic
structure, and the fundamental vector field corresponding 
to $\xi\in \Omega^0(\Sigma,\g)$ is given by 
$$ \xi_{\A(\Sigma)}(A)=-\d_A(\xi). $$ 
According to Atiyah and Bott \cite{AB,A}, a moment map $\Phi$ for this
action is given by
$$ \l\Phi(A),\xi\r=\int_\Sigma F_A\cdot\xi+ 
\int_{\partial \Sigma} \iota^*\,(A\cdot\xi),$$
where the second integral is defined with respect to the induced 
orientation on $\partial\Sigma$.  
Let $\Gpar(\Sig) \subset \G(\Sig)$ be the kernel of the restriction
map to the boundary so that there is an exact sequence
$$ 1 \ra \Gpar(\Sig) \ra \G(\Sig) \ra \G(\partial \Sig) \ra 1. $$
The moment map for the action of $\Gpar(\Sig)$ on $\A(\Sig)$
is $A \mapsto F_A$ and hence the symplectic quotient of $\A(\Sig)$ by 
$\Gpar(\Sig)$ is
$$ \M(\Sig) := \A_F(\Sig)/\Gpar(\Sig) $$
where $\A_F(\Sig) \subset \A(\Sig)$ is the space of flat connections.
If $\partial \Sig=\emptyset$ then $\M(\Sig)$ is a compact, finite
dimensional stratified symplectic space (in general singular)
\cite{SL}.  On the other hand, if $\partial \Sig\not=\emptyset$ then
according to Donaldson \cite{D} $\M(\Sig)$ is a smooth
infinite-dimensional symplectic manifold.  Just as for Riemann
surfaces without boundary, $\M(\Sigma)$ admits an alternative
description via holonomies (Theorem \ref{HolonomyDescription}).  It
has a residual Hamiltonian action of the gauge group $\G(\partial
\Sig)$ of the boundary with moment map
$$ \Phi: \M(\Sig) \ra \Omega^1(\partial \Sig,\g), \ \ [ A ] \mapsto 
\iota^* A.$$
By choosing parametrizations of the boundary components $B_i\cong S^1$
compatible with the induced orientation on $\p\Sig$ one obtains an
identification $\G(\partial \Sig) \cong (LG)^b\cong L(G^b)$ where $LG =
\text{Map}(S^1,G)$ is the loop group of $G$.  

\begin{remark}
For any complex structure on $\Sig$, the Hodge star operator gives 
rise to an $LG$-invariant complex structure 
on $\M(\Sig)$ which makes $\M(\Sig)$ into an $LG$-K\"ahler manifold.
However, in this paper we will not make use of this complex 
structure. 
\end{remark}

\begin{example}
The moduli space $\M(\Sig^1_0)$ for the 
disk  is the space $ LG /G =\Omega G$ of based
loops in $G$ (fundamental homogeneous space of $LG$).
See e.g. \cite{PS,Fr}.
\end{example}

We now explain how to construct a pre-quantum line bundle $L(\Sig)$
over $\M(\Sig)$ which carries an action of a central extension of
$\G(\partial \Sig)\cong LG^b$.
Since $\A(\Sig)$ is an affine space, the trivial line bundle
$\A(\Sig)\times \C$ with 
connection 1-form
$$ \theta_A:  \, T_A\A(\Sig)\cong\Omega^1(\Sigma,\g)\ra \R,\,\,a\mapsto 
\f{1}{2}\int_\Sig
a\stackrel{ _\cdot}{\wedge} A $$ 
is a pre-quantum line bundle over $\A(\Sig)$. 
The central $S^1$-extension $\widehat{\G(\Sig)}$  
of the gauge group defined by the cocycle 
\begin{equation}
c(g_1,g_2)=\,\exp\Big(-\f{i}{4\pi}
\int_\Sig  g_1^{-1}\d\,g_1\stackrel{ _\cdot}{\wedge}\d g_2 \,g_2^{-1}
\Big),\label{cocycle}
\end{equation}
acts on $ \A(\Sig) \times \C$ by $\theta$-preserving automorphisms via
$$ (g,z) \cdot (A,w)=(g\cdot A,\,\exp\Big(\f{i}{4\pi}\int_\Sig\,
g^{-1}\d g \stackrel{ _\cdot}{\wedge} A\Big )\,z \, w).$$
We will show in the appendix, along the lines of the papers
\cite{RSW,W2,Mi} that the extension has a canonical trivialization
over the subgroup $\G_\p(\Sig)\subset \G(\Sig)$ so that
$\G_\partial(\Sig)$ acts on $\A(\Sig) \times \C$.  The quotient
$$  L(\Sig) =(\A(\Sig) \times \C)\qu \G_\p(\Sig) = 
 (\A_F(\Sig) \times \C) / \G_\partial(\Sig) $$
with the induced connection is a pre-quantum line bundle over
$\M(\Sig)$ which carries an action of the central extension
$$ \widehat{\G(\partial \Sig)} := \widehat{\G(\Sig)}/\G_\partial(\Sig) $$
of $\G(\partial \Sig)$.  The restriction of $\widehat{\G(\partial \Sig)}$
to the boundary is isomorphic to the basic central extension
$\widehat{LG^b}$ of the loop group \cite{PS}.  We denote by
$L^\levi(\Sig)$ the $\levi$-th tensor power of $L(\Sig)$, and by 
$M^\levi(\Sig)$ the moduli space with $\levi$ times the equivariant
symplectic form, so that $L^\levi(\Sig) \ra M^\levi(\Sig)$ is an
$\widehat{LG^b}$-equivariant pre-quantum bundle.

\subsection{Gluing equals reduction}

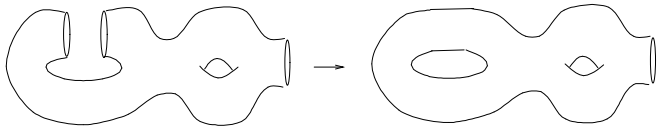
\begin{figure}[htb] 
\begin{center}
\setlength{\unitlength}{0.00033333in}
\begingroup\makeatletter\ifx\SetFigFont\undefined%
\gdef\SetFigFont#1#2#3#4#5{%
  \reset@font\fontsize{#1}{#2pt}%
  \fontfamily{#3}\fontseries{#4}\fontshape{#5}%
  \selectfont}%
\fi\endgroup%
{\renewcommand{\dashlinestretch}{30}
\begin{picture}(10144,1921)(0,-10)
\put(4391,984){\ellipse{86}{710}}
\put(948,1426){\ellipse{106}{690}}
\put(1531,1426){\ellipse{106}{734}}
\path(1564,1793)	(1613.708,1803.968)
	(1660.076,1813.618)
	(1703.284,1821.957)
	(1743.515,1828.994)
	(1815.768,1839.196)
	(1878.285,1844.294)
	(1932.516,1844.355)
	(1979.912,1839.447)
	(2021.924,1829.639)
	(2060.000,1815.000)

\path(2060,1815)	(2112.004,1778.178)
	(2156.912,1724.440)
	(2197.022,1659.397)
	(2234.631,1588.660)
	(2272.036,1517.841)
	(2311.535,1452.550)
	(2355.423,1398.400)
	(2406.000,1361.000)

\path(2406,1361)	(2455.144,1346.421)
	(2513.920,1342.321)
	(2572.486,1347.560)
	(2621.000,1361.000)

\path(2621,1361)	(2672.846,1393.590)
	(2722.453,1441.978)
	(2770.479,1500.932)
	(2817.585,1565.216)
	(2864.430,1629.599)
	(2911.675,1688.846)
	(2959.978,1737.724)
	(3010.000,1771.000)

\path(3010,1771)	(3080.017,1798.518)
	(3119.959,1810.702)
	(3162.534,1821.674)
	(3207.252,1831.315)
	(3253.621,1839.504)
	(3301.149,1846.120)
	(3349.345,1851.043)
	(3397.718,1854.151)
	(3445.776,1855.326)
	(3493.029,1854.446)
	(3538.983,1851.390)
	(3583.149,1846.038)
	(3625.035,1838.269)
	(3700.000,1815.000)

\path(3700,1815)	(3755.684,1779.168)
	(3805.904,1725.679)
	(3852.382,1660.373)
	(3896.844,1589.091)
	(3941.012,1517.673)
	(3986.612,1451.958)
	(4035.367,1397.787)
	(4089.000,1361.000)

\path(4089,1361)	(4129.611,1348.507)
	(4182.264,1344.342)
	(4253.035,1348.507)
	(4297.113,1353.712)
	(4348.000,1361.000)

\path(917,1771)	(871.591,1780.346)
	(829.236,1788.446)
	(789.768,1795.300)
	(753.022,1800.908)
	(687.033,1808.385)
	(629.941,1810.878)
	(580.419,1808.385)
	(537.140,1800.908)
	(464.000,1771.000)

\path(464,1771)	(417.533,1740.035)
	(372.318,1703.602)
	(328.604,1662.257)
	(286.639,1616.552)
	(246.672,1567.044)
	(208.951,1514.286)
	(173.726,1458.834)
	(141.244,1401.241)
	(111.754,1342.063)
	(85.506,1281.853)
	(62.747,1221.167)
	(43.727,1160.560)
	(28.693,1100.584)
	(17.895,1041.796)
	(11.581,984.750)
	(10.000,930.000)

\path(10,930)	(13.396,882.406)
	(21.789,833.178)
	(34.762,782.738)
	(51.897,731.510)
	(72.776,679.916)
	(96.981,628.378)
	(124.094,577.318)
	(153.696,527.159)
	(185.371,478.323)
	(218.700,431.232)
	(253.264,386.310)
	(288.647,343.978)
	(324.430,304.658)
	(360.195,268.774)
	(430.000,209.000)

\path(430,209)	(465.328,185.259)
	(505.526,163.107)
	(549.934,142.533)
	(597.890,123.527)
	(648.737,106.079)
	(701.812,90.177)
	(756.458,75.813)
	(812.013,62.975)
	(867.817,51.654)
	(923.211,41.838)
	(977.535,33.519)
	(1030.128,26.685)
	(1080.332,21.327)
	(1127.485,17.433)
	(1170.927,14.994)
	(1210.000,14.000)

\path(1210,14)	(1284.409,17.914)
	(1326.987,22.896)
	(1372.348,29.641)
	(1419.890,37.959)
	(1469.012,47.661)
	(1519.114,58.557)
	(1569.595,70.456)
	(1619.854,83.170)
	(1669.291,96.507)
	(1717.305,110.278)
	(1763.295,124.294)
	(1806.660,138.363)
	(1846.800,152.298)
	(1915.000,179.000)

\path(1915,179)	(1962.570,204.553)
	(2016.484,240.843)
	(2074.757,284.172)
	(2135.403,330.842)
	(2196.435,377.158)
	(2255.869,419.421)
	(2311.720,453.934)
	(2362.000,477.000)

\path(2362,477)	(2423.092,493.978)
	(2460.106,501.405)
	(2499.204,507.127)
	(2538.713,510.396)
	(2576.964,510.463)
	(2643.000,498.000)

\path(2643,498)	(2697.582,464.531)
	(2747.607,413.909)
	(2794.485,351.797)
	(2839.625,283.854)
	(2884.437,215.742)
	(2930.331,153.121)
	(2978.715,101.654)
	(3031.000,67.000)

\path(3031,67)	(3098.858,42.797)
	(3137.151,33.018)
	(3177.748,24.843)
	(3220.193,18.315)
	(3264.031,13.473)
	(3308.809,10.360)
	(3354.070,9.018)
	(3399.360,9.487)
	(3444.225,11.810)
	(3488.209,16.028)
	(3530.857,22.183)
	(3571.715,30.316)
	(3610.328,40.469)
	(3679.000,67.000)

\path(3679,67)	(3739.179,111.260)
	(3791.705,175.333)
	(3815.895,212.735)
	(3839.128,252.614)
	(3861.723,294.144)
	(3884.000,336.500)
	(3906.277,378.856)
	(3928.872,420.386)
	(3952.105,460.265)
	(3976.295,497.667)
	(4028.821,561.740)
	(4089.000,606.000)

\path(4089,606)	(4126.567,618.656)
	(4174.886,622.875)
	(4239.512,618.656)
	(4279.676,613.383)
	(4326.000,606.000)

\path(955,1064)	(911.360,1056.232)
	(871.352,1048.506)
	(801.752,1032.960)
	(745.231,1016.924)
	(700.825,999.957)
	(644.486,961.476)
	(625.000,914.000)

\path(625,914)	(655.636,852.095)
	(691.235,824.599)
	(734.954,799.816)
	(782.741,778.057)
	(830.545,759.629)
	(874.315,744.840)
	(910.000,734.000)

\path(910,734)	(971.222,720.555)
	(1044.381,711.305)
	(1084.186,708.197)
	(1125.475,706.071)
	(1167.749,704.906)
	(1210.506,704.679)
	(1253.248,705.368)
	(1295.474,706.952)
	(1336.684,709.408)
	(1376.379,712.714)
	(1449.221,721.788)
	(1510.000,734.000)

\path(1510,734)	(1546.056,743.082)
	(1590.580,754.865)
	(1639.544,769.648)
	(1688.920,787.730)
	(1734.681,809.408)
	(1772.800,834.980)
	(1810.000,899.000)

\path(1810,899)	(1785.481,967.737)
	(1719.160,1027.918)
	(1676.487,1051.392)
	(1630.759,1068.639)
	(1584.442,1078.296)
	(1540.000,1079.000)

\path(1540,1079)	(1525.000,1064.000)

\path(3025,929)	(3073.642,883.516)
	(3116.895,845.654)
	(3190.574,790.685)
	(3252.715,759.873)
	(3310.000,749.000)

\path(3310,749)	(3371.159,755.551)
	(3439.124,782.656)
	(3477.965,805.125)
	(3521.276,834.184)
	(3569.980,870.314)
	(3625.000,914.000)

\path(3100,884)	(3128.522,925.977)
	(3154.717,960.886)
	(3202.409,1011.436)
	(3247.646,1039.519)
	(3295.000,1049.000)

\path(3295,1049)	(3351.175,1041.153)
	(3406.150,1011.350)
	(3465.550,955.372)
	(3498.667,916.249)
	(3535.000,869.000)

\put(10093,1020){\ellipse{86}{710}}
\path(6625,1814)(7300,1829)
\path(6640,1169)(7150,1184)
\path(4810,914)(5260,914)
\path(5140.000,884.000)(5260.000,914.000)(5140.000,944.000)
\path(6619,1807)	(6573.912,1802.736)
	(6531.994,1798.391)
	(6457.040,1789.317)
	(6392.879,1779.488)
	(6338.256,1768.620)
	(6291.914,1756.424)
	(6252.595,1742.613)
	(6190.000,1709.000)

\path(6190,1709)	(6152.693,1681.196)
	(6113.949,1649.137)
	(6074.287,1613.203)
	(6034.228,1573.773)
	(5994.293,1531.226)
	(5955.002,1485.942)
	(5916.875,1438.301)
	(5880.432,1388.680)
	(5846.196,1337.460)
	(5814.684,1285.020)
	(5786.419,1231.739)
	(5761.921,1177.997)
	(5741.709,1124.173)
	(5726.305,1070.646)
	(5716.228,1017.795)
	(5712.000,966.000)

\path(5712,966)	(5713.792,916.944)
	(5720.996,866.631)
	(5733.149,815.445)
	(5749.793,763.768)
	(5770.465,711.986)
	(5794.705,660.480)
	(5822.053,609.634)
	(5852.048,559.832)
	(5884.228,511.458)
	(5918.134,464.894)
	(5953.304,420.523)
	(5989.279,378.730)
	(6061.796,304.410)
	(6132.000,245.000)

\path(6132,245)	(6167.328,221.259)
	(6207.526,199.107)
	(6251.934,178.533)
	(6299.890,159.527)
	(6350.737,142.079)
	(6403.812,126.177)
	(6458.458,111.813)
	(6514.012,98.975)
	(6569.817,87.654)
	(6625.211,77.838)
	(6679.535,69.519)
	(6732.128,62.685)
	(6782.332,57.327)
	(6829.485,53.433)
	(6872.927,50.994)
	(6912.000,50.000)

\path(6912,50)	(6986.409,53.914)
	(7028.987,58.896)
	(7074.348,65.641)
	(7121.890,73.959)
	(7171.012,83.661)
	(7221.114,94.557)
	(7271.595,106.456)
	(7321.854,119.170)
	(7371.291,132.507)
	(7419.305,146.278)
	(7465.295,160.294)
	(7508.660,174.363)
	(7548.800,188.298)
	(7617.000,215.000)

\path(7617,215)	(7664.570,240.553)
	(7718.484,276.843)
	(7776.757,320.172)
	(7837.403,366.842)
	(7898.435,413.158)
	(7957.869,455.421)
	(8013.720,489.934)
	(8064.000,513.000)

\path(8064,513)	(8125.092,529.978)
	(8162.106,537.405)
	(8201.204,543.127)
	(8240.713,546.396)
	(8278.964,546.463)
	(8345.000,534.000)

\path(8345,534)	(8399.582,500.531)
	(8449.607,449.909)
	(8496.485,387.797)
	(8541.625,319.854)
	(8586.437,251.742)
	(8632.331,189.121)
	(8680.715,137.654)
	(8733.000,103.000)

\path(8733,103)	(8800.858,78.797)
	(8839.151,69.018)
	(8879.748,60.843)
	(8922.193,54.315)
	(8966.031,49.473)
	(9010.809,46.360)
	(9056.070,45.018)
	(9101.360,45.487)
	(9146.225,47.810)
	(9190.209,52.028)
	(9232.857,58.183)
	(9273.715,66.316)
	(9312.328,76.469)
	(9381.000,103.000)

\path(9381,103)	(9441.179,147.260)
	(9493.705,211.333)
	(9517.895,248.735)
	(9541.128,288.614)
	(9563.723,330.144)
	(9586.000,372.500)
	(9608.277,414.856)
	(9630.872,456.386)
	(9654.105,496.265)
	(9678.295,533.667)
	(9730.821,597.740)
	(9791.000,642.000)

\path(9791,642)	(9828.567,654.656)
	(9876.886,658.875)
	(9941.512,654.656)
	(9981.676,649.383)
	(10028.000,642.000)

\path(8727,965)	(8775.642,919.516)
	(8818.895,881.654)
	(8892.574,826.685)
	(8954.715,795.873)
	(9012.000,785.000)

\path(9012,785)	(9073.159,791.551)
	(9141.124,818.656)
	(9179.965,841.125)
	(9223.276,870.184)
	(9271.980,906.314)
	(9327.000,950.000)

\path(8802,920)	(8830.522,961.977)
	(8856.717,996.886)
	(8904.409,1047.436)
	(8949.646,1075.519)
	(8997.000,1085.000)

\path(8997,1085)	(9053.175,1077.153)
	(9108.150,1047.350)
	(9167.550,991.372)
	(9200.667,952.249)
	(9237.000,905.000)

\path(7266,1829)	(7307.527,1822.299)
	(7346.127,1815.895)
	(7415.124,1803.848)
	(7474.144,1792.594)
	(7524.343,1781.870)
	(7566.873,1771.412)
	(7602.890,1760.957)
	(7660.000,1739.000)

\path(7660,1739)	(7708.841,1709.703)
	(7762.463,1667.831)
	(7819.500,1617.750)
	(7878.581,1563.826)
	(7938.340,1510.424)
	(7997.409,1461.908)
	(8054.418,1422.645)
	(8108.000,1397.000)

\path(8108,1397)	(8156.088,1386.289)
	(8214.926,1381.696)
	(8274.051,1384.756)
	(8323.000,1397.000)

\path(8323,1397)	(8374.757,1429.462)
	(8424.336,1477.758)
	(8472.380,1536.655)
	(8519.529,1600.920)
	(8566.424,1665.320)
	(8613.707,1724.623)
	(8662.019,1773.593)
	(8712.000,1807.000)

\path(8712,1807)	(8781.979,1834.407)
	(8821.897,1846.570)
	(8864.448,1857.540)
	(8909.142,1867.195)
	(8955.488,1875.408)
	(9002.996,1882.057)
	(9051.176,1887.016)
	(9099.539,1890.163)
	(9147.593,1891.373)
	(9194.849,1890.521)
	(9240.816,1887.484)
	(9285.005,1882.137)
	(9326.926,1874.357)
	(9366.087,1864.020)
	(9402.000,1851.000)

\path(9402,1851)	(9457.704,1815.276)
	(9507.932,1761.819)
	(9554.410,1696.493)
	(9598.866,1625.159)
	(9643.026,1553.680)
	(9688.618,1487.919)
	(9737.367,1433.738)
	(9791.000,1397.000)

\path(9791,1397)	(9831.765,1384.299)
	(9884.474,1380.065)
	(9955.195,1384.299)
	(9999.208,1389.591)
	(10050.000,1397.000)

\path(6655,1169)	(6610.578,1155.047)
	(6569.907,1141.504)
	(6499.333,1115.326)
	(6442.319,1089.824)
	(6397.903,1064.356)
	(6343.021,1010.959)
	(6327.000,950.000)

\path(6327,950)	(6360.712,882.565)
	(6396.439,854.837)
	(6439.688,830.900)
	(6486.663,810.600)
	(6533.570,793.785)
	(6576.614,780.303)
	(6612.000,770.000)

\path(6612,770)	(6673.222,756.555)
	(6746.381,747.305)
	(6786.186,744.197)
	(6827.475,742.071)
	(6869.749,740.906)
	(6912.506,740.679)
	(6955.248,741.368)
	(6997.474,742.952)
	(7038.684,745.408)
	(7078.379,748.714)
	(7151.221,757.788)
	(7212.000,770.000)

\path(7212,770)	(7247.932,778.906)
	(7292.155,790.265)
	(7340.747,804.523)
	(7389.788,822.129)
	(7435.354,843.531)
	(7473.526,869.179)
	(7512.000,935.000)

\path(7512,935)	(7498.808,999.738)
	(7444.265,1058.926)
	(7399.055,1088.263)
	(7340.589,1118.402)
	(7267.895,1150.071)
	(7225.909,1166.707)
	(7180.000,1184.000)

\end{picture}
}
\caption{The Riemann surfaces $\hat{\Sig}$ and $\Sig$ \label{GluingFig}}
\end{center}
\end{figure}

Let ${\Sig}$ be a compact oriented Riemann surface obtained from a
second (possibly disconnected) Riemann surface $\hat{\Sig}$ by gluing
along two boundary components $B_\pm \subset \partial \hat{\Sig}$ by
the map $ B_+\ra B_-;\,z\mapsto z^{-1}$.  (See Figure
\ref{GluingFig}.)  In this case $\M(\Sig)$ can be obtained from
$\M(\hat{\Sig})$ by a symplectic reduction.  Let $LG \ra \G(B_+)
\times \G(B_-) \subset \G(\partial \hat{\Sig})$ denote the
anti-diagonal embedding induced by the map $z \ra (z,z^{-1})$.  The
anti-diagonal embedding lifts to the central extension and therefore
acts on $L(\hat{\Sig})\ra \M(\hat{\Sig})$.  We have the following
theorem which we learned from S. Martin:
\begin{theorem} 
\label{GluingEqualsReduction} 
The moduli space $\M(\Sig)$ and line bundle $L(\Sig)$ are given by
symplectic reduction by the anti-diagonal action:
$$\M(\Sig) = \M(\hat{\Sig}) \qu LG, \ \ \ \ L(\Sig) = L(\hat{\Sig})
\qu LG.$$
\end{theorem}

On a set-theoretical level, this follows easily by noting that 
the moment map for the anti-diagonal action is given by 
$$ [\hat{A}] \mapsto \iota_+^* \hat{A} -\iota_-^* \hat{A} $$
where $\iota_\pm:\,S^1 \ra \hat{\Sig}$ is the composition of the map
$z \mapsto z^{\pm 1}$ with the inclusion $B_\pm \ra \hat{\Sig}$.
Therefore, if an equivalence class $[\hat{A}]$ is in the zero level
set then the pullbacks of the representative $\hat{A}$ to $B_\pm$ are
equal.  In this case one can choose a representative $\hat{A}$ so that
the two sides patch together to form a connection $A$ on $\Sigma$. We
give a rigorous argument in Appendix \ref{theproof}.

\subsection{Factorization}

Let $G$ be a compact Lie group and $M$ a Hamiltonian $G$-space.
Suppose that by some quantization procedure one can construct out of
these data a virtual representation $Q(M)$ of $G$.  The ``quantization
commutes with reduction'' principle (as formulated by
Guillemin-Sternberg \cite{GS2}) says that in this case the
quantization of the symplectic reduction $M \qu G$ should equal
the invariant part of the quantization:
\begin{equation} \label{Commutes}
 Q(M\qu G) = Q(M)^G. 
\end{equation}
In the context of K\"ahler quantization (\ref{Commutes}) was proved in
\cite{GS2} for smooth symplectic quotients and in \cite{Sja} for
singular quotients.  Another approach (which does not require the
existence of an invariant K\"ahler structure) is to define $Q(M)$ to
be the equivariant index (Riemann-Roch number)
of the $\text{Spin}_c$-Dirac
operator $\dirac_L$ associated to a compatible, invariant almost
complex structure on $M$ and pre-quantum line bundle $L \ra M$ (see
Section \ref{RiemannRochNumbers})
$$ Q(M) = \RR(M,L).$$
In this form, the principle has been proved in
\cite{Gu,M1,Vg1} in the abelian case and in \cite{M2} for the non-abelian 
case; the case of singular quotients has been dealt with in \cite{MS}.

Heuristically, the factorization property for moduli spaces of flat 
connections follows from an application of
(\ref{Commutes}) in the setting of Hamiltonian actions of loop groups.
Let $\hat{\Sig}$ be a compact oriented Riemann surface (possibly
disconnected) and $\Sig$ the Riemann surface formed by gluing along
two boundary components $B_\pm \subset \partial \hat{\Sig}$.  Theorem
\ref{GluingEqualsReduction} and (\ref{Commutes}) would imply an
isomorphism of $\widehat{\G(\partial\Sig)}$-representations
$$ Q(\M^\levi(\Sig)) = Q(\M^\levi(\hat{\Sig}))^{LG} \ \ \text{(Factorization)}.
$$
We want to emphasize that we do not prove the factorization theorem in
this form, which would require the construction of the quantization of
an infinite dimensional symplectic manifold.

Rather, note the following two corollaries of the principle for
compact groups $G$.  For any dominant weight $\mu$, let $*\mu$ be the
dominant weight for the dual representation $V_{\mu}^*$, which by
Borel-Weil can be realized as the quantization of the coadjoint orbit
$\O_{* \mu} = G \cdot *\mu$.  The reduction $M_\mu := M \times \O_{*
\mu} \qu G$ is called the reduction of $M$ at $\mu$.  As a corollary of
the principle (\ref{Commutes}) one has that $ Q(M_\mu) = (Q(M)\otimes
V_\mu^*)^G$, so that 
$$ Q(M) = \oplus_\mu \, Q(M_\mu) \, V_\mu .$$
Now suppose that $M$ is a compact quantizable Hamiltonian $G \times
G$-space, and let $G$ act on $M$ by the diagonal action.  Then
$$ Q(M \qu  G) = Q(M)^G = \bigoplus_{\mu,\nu} Q(M_{\mu,\nu})\Big(  V_\mu \otimes
V_\nu \Big)^G = \bigoplus_\mu Q(M_{\mu,*\mu}).
$$
Our main result is a generalization of this formula to the setting of
Hamiltonian $LG$-actions with proper moment maps:

\begin{theorem} \label{MainResult} 
Let $G$ be a compact connected simply-connected Lie group, and $\Alc$
the corresponding fundamental alcove.  Let $M$ be a Hamiltonian 
$L(G\times G)$-Banach manifold with proper moment map at non-zero level
$\levi$ and $\widehat{LG^2}$-equivariant
pre-quantum line bundle $L \ra M$. Then
$$ \RR(M \qu LG,L \qu LG)=\sum_{\mu\in m \Alc \cap
 \Lambda^*}\, \RR(M_{\mu,*\mu},L_{\mu,*\mu}). $$
\end{theorem}
The properness assumption guarantees that all quotients are finite
dimensional and compact. Their Riemann-Roch numbers can be defined using
desingularization if necessary (see Section \ref{DesingSec}). 

As a special case Theorem \ref{MainResult} gives
\begin{theorem}[Factorization] \label{FusionRules} 
Let $\hat{\Sig}$ be a compact oriented 
Riemann surface (possibly disconnected) with $b \ge 2$ boundary
components and $\Sig$ the Riemann surface formed by gluing along two
boundary components $B_\pm \subset \partial \hat{\Sig}$.  Given a
level $\levi \in \N$ and dominant weights $\nu =
(\nu_1,\ldots,\nu_{b-2})$ at level $\levi$, one has
$$ \RR(\M^{\levi} (\Sig)_\nu,L^{\levi}(\Sig)_\nu)=\sum_{\mu\in
 \levi \Alc \cap \Lambda^*}\, \RR(\M^{\levi}(\hat{\Sig})_{\nu,\mu,*
 \mu},L^{\levi}(\hat{\Sig})_{\nu,\mu,* \mu}) .$$
\end{theorem}

A gauge-theoretic interpretation of the reductions $\M^m(\Sig)_\nu$ is
given as follows.  The fundamental alcove $\Alc$ can be viewed as the
set of conjugacy classes $\Alc\cong G/\Ad(G)$. The Cartan
involution $*:\Alc\to\Alc$ sends the conjugacy class $(g)$ to
$(g^{-1})$.  For any Riemann surface $\Sig$ with $b$ boundary
components, the reduced space $\M^m(\Sig)_\nu$ for
$\nu=(\nu_1,\ldots,\nu_b)\in m\,\Alc^b$ is the moduli space of flat
connections on the {\bf marked Riemann surface}
$(\Sig;\nu_1,\ldots,\nu_b)$, for which the holonomies $g_j\in G$
around the boundary components $B_j$ are required to lie in the
conjugacy classes $\nu_j/m$.

As examples we review two well-known applications of the factorization theorem.
\subsection{Fusion product}
Let $\Sig$ be the three-holed sphere and $\M(\Sig)$ the associated
$LG^3$-Hamiltonian manifold.  For $\mu,\nu,\alpha \in
\levi \Alc\cap\Lambda^*$  define 
\begin{equation}\label{FusionCoefficients}
N^{\levi}_{\mu,\nu;\alpha}:= \RR(\M^\levi(\Sig)_{\mu,\nu,*\alpha},
L^\levi(\Sig)_{\mu,\nu,*\alpha}) 
\end{equation}
the Riemann-Roch number of the moduli space corresponding to the
markings $\mu,\nu,*\alpha$.  The integers $N^{\levi}_{\mu,\nu,\alpha}$
are a set of {\bf fusion rules}. This means that they satisfy the
axioms
$$ N^\levi_{\mu,\nu;\alpha}=N^\levi_{\nu,\mu;\alpha}=
N^\levi_{*\mu,*\nu;*\alpha},$$
$$ N^\levi_{\mu,*\nu;0}=N^\levi_{\mu,0;\nu}=\delta_{\mu,\nu}, $$
and   
\begin{equation}\label{Associativity}
 \sum_{\alpha} N^\levi_{\mu,\alpha;\nu}\,N^\levi_{\,\rho;\alpha}
= \sum_{\alpha} N^\levi_{\beta,\alpha;\mu}\,N^\levi_{\nu,\rho;\alpha}.
\end{equation} 
That (\ref{FusionCoefficients}) satisfies the 
associativity property (\ref{Associativity}) follows because by
Theorem \ref{FusionRules}, both sides are equal to
$$
\RR(\M^\levi(\ti{\Sig})_{*\mu,\nu,\beta,\rho},
L^\levi(\ti{\Sig})_{*\mu,\nu,\beta,\rho})
$$ 
where $\ti{\Sig}$ is the four-holed sphere.

The {\bf fusion product} on the free group $\Rep_m(LG)$ of 
irreducible $LG$-representations at level $m\in\N$
is the commutative associative product 
$$ V_\mu \fus V_\nu = \bigoplus_\alpha N^{\levi}_{\mu,\nu;\alpha}
V_{\alpha}. $$
The holomorphic induction map at level $m$
$$\Ind_m: (\Rep(G),\otimes) \ra (\Rep_m(LG),\fus)$$
is a homomorphism of fusion rings.  This statement is proved for
example in Teleman \cite{Tel} and Faltings \cite{F}.  In particular,
the fusion coefficients at level $\levi$ can be expressed in terms of
the coefficients for the tensor product on $\Rep(G)$.

\subsection{Verlinde formula}
As noted by Verlinde, in case $b = 0$ and $g \ge 2$ factorization
leads to an expression for the quantum dimension of $\M(\Sig)$ in
terms of eigenvalues of a certain symmetric matrix.  In this case,
$\Sig$ can be obtained by gluing together $g-1$ copies of the
two-punctured torus $\Sig^2_1$.  Let $N = \# \Alc \cap \levi
\Lambda^*$ be the number of points in the fundamental alcove and $A$
the $N \times N$ matrix with coefficients
$$ A^\alpha_\beta = \RR(\M^{\levi}
(\Sig^2_1)_{\alpha,*\beta},L^{\levi}(\Sig^2_1)_{\alpha,*\beta}) =
\sum_{\mu,\nu} N^{\levi}_{\mu,\nu;\alpha}
N^{\levi}_{\mu,\nu;\beta}.$$
Theorem \ref{MainResult} implies that 
$$  \RR(\M^{\levi} (\Sig)_\nu,L^{\levi}(\Sig)_\nu)
= \text{Tr}(A^{g-1}) = \sum \lambda_i^{g-1} $$
where $\lambda_i$ are the eigenvalues of $A$.  For information on how
to obtain the explicit Verlinde formula from this approach, see
Beauville \cite{B}.

\section{Hamiltonian loop group actions}

Let $G$ be a connected, simply connected compact Lie group,  
$T$ a maximal torus and $W=N_G(T)/T$ the Weyl group. We normalize the 
invariant inner product on $\g$ as explained in Section 
\ref{ConstructionModuliSpace}, thereby identifying $\g\cong\g^*$ 
and $\t\cong\t^*$. The integral lattice $\exp^{-1}(1)\subset \t$ 
will be denoted by $\Lambda$, the weight lattice by $\Lambda^*$, 
and the affine Weyl group by $W_{\text{aff}}=W \ltimes \Lambda$. 

Let the loop group $LG$ 
be the Banach Lie group consisting of maps $S^1\ra G$ of some fixed
Sobolev class $s>\f{1}{2}$.  (By the Sobolev embedding theorem, $LG$
consists of continuous loops, so that multiplication and inversion are
defined pointwise.)  
Recall \cite{PS} that the Lie algebra of the
central extension $\widehat{LG}$ is the product
$\widehat{L\g}=L\g\times\R$, with bracket
$$ [(\xi_1,t_1),(\xi_2,t_2)]=
\Big([\xi_1,\xi_2],\,\f{1}{2\pi}\oint \xi_1\cdot\xi_2'\Big)$$
where the prime indicates the derivative with respect to the $S^1$
coordinate.  We define $L \g^*$ to be the space of maps from $S^1$ to
$\g^*$ of Sobolev class $s-1$.  The natural pairing of $L \g^*$ with
$L \g$ given by integration makes $L \g^*$ into a subset of the
topological dual of $L \g$.  The coadjoint action of $LG$ on
$\widehat{L\g^*}:=L\g^*\times\R$ is given by
\begin{equation} \label{CoadjointAction}
 g\cdot (\xi,\lambda)=(\Ad_g(\xi)-\lambda\,g'\,g^{-1},\lambda).
\end{equation}
It follows that under the identification of $L\g^*$ with $\g$-valued
$1$-forms of Sobolev class $s-1$ given by the Hodge $*$-operator and
the inner product on $\g$, the $LG$-actions on elements of $L\g^*$
considered as connections correspond to the action on the affine
hyper-plane $L\g^* \times \{1\}\subset \widehat{L\g^*}$.  Henceforth
we fix a level $\lev\not=0$ and identify $L\g^*$ with the affine
hyper-plane $L\g^* \times \{\lev \}$.  We denote by
$$ \Hol:\,L\g^*\ra G$$
the map that sends $\xi\in L\g^*$ to the holonomy of $\xi/\lev$,
considered as a connection 1-form on $S^1$. Then $\Hol$ has the
equivariance property $\Hol(g\cdot\xi)=\Ad_{g(0)}\Hol(\xi)$. The
restriction of $\Hol$ to $\g \cong \g^* \subset L\g^*$ is given by
$\Hol(\xi)=\exp(2\pi\xi/\lev)$.

We choose a closed positive Weyl chamber
$\t_+\subset\t$ and let 
$\Alc\subset \t_+$ be the corresponding fundamental alcove. 
There are natural identifications 
$$ \t_+\cong \t/W \cong \g/\Ad(G),$$ 
that is every (co)adjoint orbit meets the positive Weyl chamber 
in exactly one point. Similarly, for the  affine $LG$-action on $L\g^*$
at level $\lev$ and the 
action of 
$W_{\text{aff}}=W\ltimes \Lambda$  on $\t$  given by 
$(w,v)\cdot \xi=w\cdot\xi +\lev\, v$ one 
has 
$$   \lev\Alc \cong \t/W_{\text{aff}} \cong L\g^*/LG,$$ 
that is, every coadjoint $LG$-orbit at level $\lev$ meets 
$\lev\Alc\subset \t\subset L\g^*$ in exactly one point.

By a symplectic structure on a Banach manifold $M$ we mean a closed
two-form $\om$ that is weakly non-degenerate, that is, for any $m \in
M$ the map $T_mM \ra T_m^*M$ induced by $\om$ is injective (see
e.g. \cite{AMR}).  A smooth action of $LG$ on $M$ which preserves the
symplectic form will be called Hamiltonian at level $\lev \in\R$ if
there exists a moment map $\Phi:M \ra L \g^*$ such that the
composition of $\Phi$ with the inclusion $ L \g^* \ra \widehat{L
\g^*}$ at level $\lev$ is $LG$-equivariant.  Such an action can be
considered an action of $\widehat{LG}$ with the central circle acting
trivially with constant moment map $\lev$.  We emphasize that we
require the moment maps for the $LG$-actions to be ``sufficiently
smooth'', i.e. to take values in $L\g^*$, which consists of loops of
Sobolev class $s-1$.

A $\widehat{LG}$-equivariant line bundle $L \ra M$ with connection is
{\bf pre-quantum} if its curvature is equal to the symplectic form and
the action satisfies the pre-quantum condition (\ref{LiftFormula}).
This requires in particular that $\lev$ is an integer.

\begin{example}\label{CoadjointOrbits}
The basic examples of Hamiltonian $LG$-spaces (at level $\lev\not=0$)
are coadjoint orbits $\O_\xi=LG\cdot(\xi,\lev)$.  Letting
$\delta_{\xi/\lev}$ be the covariant derivative with respect to the
connection $\f{1}{\lev}\xi$,
\begin{equation}\label{MapDelta}
\delta_{\xi/\lev}:\,L\g\ra L\g^*,\,\,\eta\mapsto \eta'+\f{1}{\lev}[\xi,\eta], 
\end{equation}
the fundamental vector field $\eta_{L\g^*}$ 
for the infinitesimal action of $L\g$ on $L\g^*$ at level $\lev$ is given by 
$\eta_{L\g^*}(\xi)=-\lev\delta_{\xi/\lev}(\eta)$. 
The symplectic form $\nu_\xi$ on $\O_\xi$ is given by the usual 
KKS formula 
\begin{equation} \label{SympForm}
\nu_\xi((\eta_1)_{L\g^*},\,(\eta_2)_{L\g^*})
=\f{\lev}{2\pi}\oint \eta_1\cdot \delta_{\xi/\lev}(\eta_2) \end{equation}
with moment map as usual the inclusion into $L\g^*$. 

For $\xi\in\Lambda^*$, the orbit $\O_\xi,\, \xi \in \lev \Alc$ admits a 
pre-quantum line bundle
$\Xi(\O_\xi)$ if and only if $\lambda=m\in \Z$ and $\xi\in \Lambda^*$. 
If $m\in\N$, 
geometric quantization of $\Xi(\O_\xi)\ra
\O_\xi$ by the Borel-Weil construction \cite{PS} gives the irreducible
positive energy representation of $LG$ at level $m$ with highest
weight $\xi$.
\end{example}

\begin{remark}\label{Inversion}  
\begin{enumerate}
\item
The inversion map $I^*: \, LG \ra LG$ transforms 
a Hamiltonian action with moment map $(\Phi,\lambda)$ into one with moment
map $(\Phi,-\lambda)$.  Hence if $M$ is a Hamiltonian $LG \times LG$-manifold
with moment map $(\Phi_+,\lambda;\Phi_-,\lambda)$ 
at level $\lambda\not=0$ then the 
anti-diagonal action of $LG$ is at level $0$, 
with moment map $\Phi_+ +\Phi_-$.
\item
\label{Rescaling}
If $(M,\omega,\Phi)$ is a Hamiltonian $LG$-Banach manifold 
at non-zero level $\lev\not=0$, then $(M,\lev^{-1}\omega,\lev^{-1}\Phi)$
is a Hamiltonian $LG$-space at level $+1$. Henceforth, we will always 
take the level to be $+1$ unless specified otherwise, and identify $L\g^*$ 
with the affine hyper-plane $L\g^*\times\{1\}\subset\widehat{L\g^*}$. 
\end{enumerate}
\end{remark}

\begin{example}\label{ModuliTwoPuncturedSphere}
Let $G$ be connected and simply connected.  Let $T^*\widehat{LG}$ be
the cotangent bundle of the central extension of $LG$. Trivialization
by left-invariant one-forms gives a diffeomorphism
$T^*\widehat{LG}\cong \widehat{LG}\times \text{Hom}(\widehat{L\g},\R)$
where $\text{Hom} (\widehat{L\g},\R)$ is the topological dual of
$L\g$. The subset $\hat{X}=\widehat{LG}\times \widehat{L\g^*}$ is a
Hamiltonian $\widehat{LG}\times \widehat{LG}$-space, with actions given
by
$$  L_{\hat{a}} (\hat{g};\xi,\lambda)=(\hat{a}\,\hat{g};\xi,\lambda),\ \ \
R_{\hat{a}}(\hat{g};\xi,\lambda)=(\hat{g}\,\hat{a}^{-1};\,
\Ad_{\hat{a}}(\xi,\lambda))$$
and moment maps 
$$\hat{\Phi}^{(L)}(\hat{g};\xi,\lambda)=\Ad_{\hat{g}}(\xi,\lambda),\ \ \
\hat{\Phi}^{(R)}(\hat{g};\xi,\lambda)=-(\xi,\lambda).$$
Let $X$ be the reduction of $\hat{X}$ 
with respect to central 
circle $S^1$ at moment level $1$. Then $X\cong LG\times L\g^*$,
and the induced left and right actions of $LG$ 
$$ L_a(g,\xi)=(a\,g,\,\xi),\ \ \ R_{a}(g,\xi)=(g\,a^{-1},\,a\cdot\xi) $$
are Hamiltonian, with moment maps
$$\Phi^{(L)}(g,\xi)=g\cdot\xi,\ \ \ \Phi^{(R)}(g,\xi)=-\xi.$$
(Here we can use the involution $I^*$ to make $\Phi^{(R)}$ into 
a moment map at level $+1$.)
It has the property that for every Hamiltonian $LG$-Banach 
manifold $M$, the reduced space $M\times X\qu LG$ by the 
anti-diagonal action is symplectomorphic to $M$ itself.

The trivial line bundle $L_{\hat{X}}=\hat{X}\times\C$ is a pre-quantum
line bundle for $\hat{X}$, with the restriction of the canonical
1-form on $T^*\widehat{LG}$ defining a pre-quantum connection;
reduction with respect to the $S^1$-action gives an $L(G\times
G)$-equivariant pre-quantum bundle $L_X\ra X$ isomorphic to the
pull-back of $\widehat{LG}\times_{S^1}\C$. If $L\ra M$ is an
$LG$-equivariant pre-quantum line bundle, the anti-diagonal reduction
$L\boxtimes L_X\qu LG$ is isomorphic to $L$ itself.

We shall see later (Example \ref{ModuliTwoPuncturedSphere2}) that $X$
is equivariantly symplectomorphic to the moduli space $\M(\Sig^2_0)$
of the annulus and that $L_X$ is equivariantly isomorphic to
$L(\Sig^2_0)$.
\end{example}

Hamiltonian loop group actions on Banach manifolds with proper moment
maps at positive level behave in many respects like Hamiltonian
actions of compact groups on finite dimensional symplectic manifolds.
This is due to the existence of finite-dimensional ``symplectic
cross-sections''.  These are symplectic analogs of highest-weight
modules in representation theory.

\subsection{Cross-sections for Hamiltonian actions of compact groups}
\label{CrossSections}
First we review the symplectic cross-section theorem in the setting of
compact, connected Lie groups $G$. Choose a maximal torus $T \subset
G$ and positive Weyl chamber $\g^*/G=\t^*_+ \subset \t^*$.

For every open face $\sigma$ of $\t^*_+$, the stabilizer subgroup
$G_\xi\subset G$ does not depend on the choice of $\xi\in\sigma$, and
is denoted by $G_\sigma$.  Since $G_\sig$ contains the maximal torus
$T$, there is a unique $G_\sigma$-invariant splitting of Lie algebras
$$ \g=[\g_\sigma,\g_\sigma]\oplus \z(\g_\sigma)\oplus \g_\sigma^\perp $$
where $\z(\g_\sigma)$ is the center of $\g_\sigma$ and 
$\g_\sigma^\perp$ is a $G_\sig$-invariant complement. 
In fact, using an invariant inner product 
to identify $\g^*$ with $\g$, $\g_\sig$ is characterized as the kernel of 
the map  $\text{ad}_\xi=[\xi,\cdot]$, and the tangent space at $\xi$ to 
the coadjoint orbit with its image; moreover $\z(\g_\sig)$ gets identified 
with the tangent space to $\sig$.
  
Let $U_\sigma\subset\g^*$ be the $G_\sigma$-invariant open subset
of $\g_\sigma^*$ defined by
$$ U_\sigma:=G_\sigma\cdot\bigcup_{\tau\subset\ol{\sigma}}\tau. $$
Then $G\cdot U_\sigma\cong G\times_{G_\sigma} U_\sigma$, which implies
that $U_\sigma$ is a slice for the coadjoint action at any
$x\in\sig$.  Now let $N$ be a Hamiltonian $G_\sigma$-space with moment
map $\Phi_N$.  The {\bf symplectic induction} construction \cite{GS1}
shows that if $\Phi_N(N)$ is contained in $U_\sigma$ then there exists
a unique symplectic structure on the associated bundle
$$\text{Ind}_{G_\sigma}^G(N):=G\times_{G_\sigma}N $$
such that the $G$-action is Hamiltonian and the symplectic form $\om$
and moment map $\Phi$ restrict to the given symplectic form and moment
map on $N$.  Given a second group $K$ with a Hamiltonian action on $N$
such that this action commutes with $G_\sigma$, one obtains a
Hamiltonian $G\times K$-action on $\text{Ind}_{G_\sigma}^G(N)$. In
particular, the action of $G$ always extends to an action of $G\times
Z(G_\sig)$, where $Z(G_\sig)$ is the center of $ G_\sig$.  The moment
map for the $Z(G_\sigma)$-action is given by the composition of $\Phi$
with the quotient map $q:\,\g^*\ra\t^*_+$, followed by projection $\t^*\ra
\z(\g_\sig)^*$.
\label{ToricSection}

The {\bf symplectic cross-section theorem} \cite{GS1} asserts that
conversely, for every Hamiltonian $G$-manifold $M$ with moment map
$\Phi:\,M \ra \g^*$, the pre-image $Y_{\sigma}=\Phinv( U_\sig)$ is a
$G_\sig$-invariant symplectic submanifold.  Consequently, the action
of $G$ on $G\cdot Y_{\sigma}$ extends to a Hamiltonian action
of $G\times Z(G_\sigma)$.  We refer to
these actions (due to Guillemin-Sternberg \cite{GS3}) as the 
{\bf induced (toric) actions} and to the map 
\begin{equation}\label{Toric1}
\ti{\Phi}=q\circ \Phi:\,M\ra \t^*
\end{equation}
as the {\bf induced (toric) moment map}.  Notice that if $E \ra M$ is
a $G$-equivariant vector bundle, the action of $Z(G_\sigma)$ on the
restriction $E \vert \Phi^{-1}( U_\sigma)$ extends by $G$-equivariance
to an action of $G\times Z(G_\sigma)$ on the restriction of $E$ to
$G\cdot\Phi^{-1}( U_\sigma)$.
\label{ToricLift}

Suppose now that $M$ is a Hamiltonian $G \times G$-manifold (e.g. a
product of Hamiltonian $G$-manifolds) with moment map $\Phi =
(\Phi_+,\Phi_-)$. Let $M_0=M\qu G$ be the symplectic quotient by the
diagonal action, and $\ti{\Phi}_0:\,M_0\to \t^*_+$ the residual toric
moment map induced by $\ti{\Phi}_+$.  The symplectic cross-section
$Y_{\sig,-\sig}=\Phinv(U_\sigma \times - U_\sigma)$ is a Hamiltonian
$G_\sig\times G_\sig$-space, with moment map the restriction of
$\Phi$. There is canonical isomorphism of (possibly singular)
symplectic quotients
$$ M_0\supset \ti{\Phi}_0^{-1}(U_\sig)\cong Y_{\sig,-\sig}\qu G_\sig.$$

\subsection{Cross-sections for Hamiltonian actions of loop groups}

We now turn to the discussion of symplectic induction and
cross-sections for Hamiltonian $LG$-actions on symplectic Banach
manifolds, where $G$ is connected and simply connected.  
We first need to summarize some  properties of the coadjoint 
action of $LG$ on $L\g^*$ at level $1$. Recall \cite{PS} 
that the evaluation map $LG\to G,\,g\mapsto g(0)$ maps the 
isotropy group $(LG)_\xi$ of a point $\xi\in L\g^*$ isomorphically 
to the stabilizer $G_{\Hol(\xi)}$ of the holonomy of $\xi$, 
in particular $(LG)_\xi$ is compact and connected
\footnote{The stabilizer groups for the adjoint
action are connected by a result of Bott-Samelson described in
\cite{BT}.}.  
For points $\xi\in\Alc\subset \t\subset L\g^*$, the inverse map is 
given by $G_{\Hol(\xi)}\mapsto (LG)_\xi,\,k\mapsto
\Ad_{\exp(-\theta\xi)}\,k$. 

It follows that the isotropy group $(LG)_\xi$ of a point $\xi\in\Alc$
contains $T\subset LG$, and depends only on the open face
$\sigma\subset \Alc$ containing it.  We denote this group by
$(LG)_\sigma$, the group $G_{\Hol(\xi)}$ by $\ol{(LG)_\sigma}$ and the
restriction of the central extension $\widehat{LG}$ to $(LG)_\sig$ by
$\widehat{(LG)}_\sig$. Note that  
$$\tau\subset\ol{\sigma} \Rightarrow (LG)_\sigma\subset(LG)_\tau. $$ 

If $\ol{\sigma}$ contains $0$ then
$(LG)_\sigma$ is contained in the subgroup $G \subset LG$ of constant
loops and is equal to the stabilizer group of points in $\sig$ 
under the coadjoint action of $G$ on $\g^*$.

Another consequence of the above description is that $(LG)_\sig$
consists only of smooth maps. Alternatively, this follows from the fact 
that its Lie algebra $(L \g)_\sig$
of $(LG)_\sig$ is equal to the kernel of the elliptic operator
$\del_\xi:\,L\g\ra L\g^*$ defined in (\ref{MapDelta}).  The image $\im
( \delta_\xi)$ is equal to the tangent space to the $LG$-orbit through
$\xi$, i.e. to the annihilator $(L \g)_\sig^0$, so that there are
$(LG)_\sig$-invariant direct sum (Hodge) decompositions
\begin{equation}
L\g^*=(L \g)_\sig^*\oplus (L \g)_\sig^0 
\label{LgSplitting}
\end{equation}
and 
\begin{equation}
L\g=(L \g)_\sig\oplus (L \g)_\sig^\perp.
\label{LieSplitting}
\end{equation}
The map $\delta_\xi$ induces a $(LG)_\sig$-equivariant 
Banach space isomorphism $(L \g)_\sig^\perp\cong(L \g)_\sig^0$.  

\begin{remark}
It is important to note that since the $LG$-action on $L\g^*$ is only
affine-linear, the action of $(LG)_\sig$ on
$(L \g)_\sig^*\cong(L \g)_\sig^*\times \{ 1 \} \subset \widehat{L\g^*}$ in
this splitting is not the coadjoint action unless $0\in\ol{\sig}$.
Indeed, from the description of the isomorphism $(LG)_\sig \cong
G_{\Hol(\xi)}$ and the coadjoint action of $\widehat{LG}$ given in
(\ref{CoadjointAction}) one finds that the action is given by
\begin{equation} \label{AffineAction}
(LG)_\sig \times (L \g)_\sig^*\ra (L \g)_\sig^*,\ \ 
(k,\eta)\mapsto (\Ad_{k^{-1}})^* (\eta - \mu) + \mu, 
\end{equation}
where $\mu$ is any element in the 
affine span of $\sig$. 
\end{remark}

As above, we now define, for every open face $\sigma\subset\Alc$
$$ U_\sigma=(LG)_\sigma\cdot \bigcup_{\tau\subset\ol{\sigma}}\tau.$$
Note that $U_\sigma$ is an open subset of $(L \g)_\sigma^*$, in particular 
it is finite dimensional and consists only of smooth elements of 
$L\g^*$. 
\begin{lemma}
The set $U_\sigma$ is a slice for all $\xi\in\sigma$ 
for the action of $LG$, i.e. the canonical map 
$$LG\times_{(LG)_\sigma} U_\sigma\ra LG\cdot U_\sigma$$
is a diffeomorphism of Banach manifolds.
\end{lemma}

\begin{proof}
The map is bijective because for any $\eta \in U_\sig$, the stabilizer
$(LG)_\eta\subset (LG)_\sigma$.
That the differential is an isomorphism
follows from the splitting (\ref{LgSplitting}).
\end{proof}

\begin{remark}\label{CentralElements}
Let $\sig = \{ \xi \} \subset \Alc$ be a vertex such that $\exp(-2\pi
\xi)$ is contained in the center $Z(G)$ of $G$, i.e.
$\ol{(LG)_\sigma}=G$. Then $f_\xi(\theta):=\exp(-\theta\xi)$ defines an
exterior automorphism of $LG$, 
$$ (f_\xi\cdot g)(\theta)= 
\Ad_{f_\xi(\theta)}\,g(\theta).$$  
Similarly, there is an automorphism of $L\g$ by 
$$
{f}_\xi\cdot \eta =\Ad_{f_\xi(\theta)}(\eta)+\xi.$$ 
These two automorphisms are compatible, that is 
$${f}_\xi\cdot(g\cdot\eta)=(f_\xi\cdot g)\cdot\,({f}_\xi\cdot\eta).$$
It follows that the slice $U_\sigma$ is isomorphic to 
the slice $U_{\{0\}}$. For example, if $G=SU(n)$ all vertices of 
$\Alc$ exponentiate to elements of the center.
\end{remark}\vskip .1in
\begin{lemma} \label{SplitLemma}
For all $\xi \in U_\sig$, the tangent space 
$T_\xi(LG\cdot\xi)\cong\im  (\delta_\xi)$ to the coadjoint orbit
$\O_\xi$ 
through $\xi$ decomposes into an $\nu_\xi$-orthogonal direct sum of 
closed symplectic subspaces,
$$ T_\xi(LG\cdot\xi)\cong T_\xi((LG)_\sig\cdot\xi)\oplus (L \g)_\sig^0, $$ 
where $(L \g)_\sig^0$ is the annihilator of $(L \g)_\sig$ in $L\g^*$.
\end{lemma}

\begin{proof}
It follows easily from the definition (\ref{SympForm}) of the
symplectic form on $LG \cdot \xi$ that the two subspaces are
symplectically orthogonal.
The proof is completed by noting that the coadjoint orbit $(LG)_\sig
\cdot \xi $ with the KKS form is a symplectic
submanifold.
\end{proof}

\begin{theorem} (Symplectic induction)
Let $\sigma$ be a face of $\Alc$ and $N$ a symplectic Banach
manifold with a Hamiltonian $(LG)_\sigma$-action and moment map
$\Phi_N:\,N\ra (L \g)_\sigma^*$ such that for some $\mu$ in the affine 
span of  $\sig$, the
image $(\Phi_N + \mu)(N)\subset U_\sigma$.  Then there exists a
unique $LG$-invariant symplectic form $\omega$ on the Banach manifold
$$\Ind_{(LG)_\sigma}^{LG}(N):=LG\times_{(LG)_\sigma}N $$
such the $LG$-action is Hamiltonian, with moment map
$\Phi:\,\Ind_{(LG)_\sigma}^{LG}(N)\ra L\g^*$, and such that the pullback
of $\omega$ (resp. $\Phi$) to $N$ is equal to the given symplectic
form (resp. moment map $\Phi_N + \mu$) on $N$.
\end{theorem}

\begin{proof}
Let $M:= {LG}\times_{(LG)_\sig} N$.  By Equation (\ref{AffineAction}) the
map $\Phi_N + \mu:\,N\ra (L \g)_\sig^*\oplus\{1\}$ extends to a unique
$LG$-equivariant map $\Phi:\,M \ra LG\times_{(LG)_\sig} (L \g)_\sig^* \hra
L\g^*$.  For $\Phi$ to be a moment map for $\omega$, one must have
\begin{equation}
\omega_x(\eta_M,\zeta_M)=\nu_{\Phi(x)}(\eta_{L\g^*},\zeta_{L\g^*}).
\label{Condition}
\end{equation}
This condition, together with $LG$-invariance of $\omega$ and the
condition that the pull-back to $N$ be $\omega_N$ completely
determines $\omega$, and also implies that $\omega$ is closed.
To show $M$ is symplectic let $x\in N\subset M$, and
$\tau\subset\Alc$ the open face containing $\Phi_N(x)$.  
By (\ref{Condition}) together with Lemma \ref{SplitLemma},
there is a
natural $\omega$-orthogonal splitting $T_x M\cong T_x N\oplus
(L \g)_\sig^0$.
Since $\sigma\subset\ol{\tau}$, this
is in fact a symplectic subspace, which shows that $T_x M$ is
symplectic.
\end{proof}

\begin{theorem}(Symplectic cross-section) \label{CrossSectionThm}
Let $(M,\om)$ be a symplectic Banach manifold, and $LG\times M\ra M$ a
Hamiltonian $LG$-action with moment map $\Phi:\,M\ra L\g^*$. For every
open face $\sigma\subset \Alc$, the {\bf symplectic cross-section}
$Y_{\sig} := \Phi^{-1}(U_{\sigma})$ is a symplectic
$(LG)_\sigma$-invariant Banach submanifold, and the action of $(LG)_\sigma$
is Hamiltonian.  The restriction $\Phi|Y_\sig$ is a moment map for the
action of $\widehat{(LG)_\sig}$, and a moment map for the $(LG)_\sig$-action
is given by $\Phi|Y_\sig- \mu$, for any $\mu$ in the affine span
of $\sig$.
\end{theorem}

\begin{proof}
By equivariance and since $U_\sigma$ is a slice for the $LG$-action,
$\Phi$ is transversal to $U_\sigma$. The implicit function theorem for
Banach manifolds thus shows that $N:=\Phi^{-1}(U_\sig)$ is a smooth
Banach submanifold.  Since $T_{\Phi(x)}\,U_\sig =(L \g)_\sig^*\subset
L\g^*$, the tangent space $T_x N$ at some $x \in N$ is equal to
$(\d_x\Phi)^{-1} ((L \g)_\sig^*)$.  In order to show that $T_xN$ is
symplectic, we show that $T_x M \cong E\oplus T_xN$ where $E$ is a
closed symplectic complement to $T_x N$ that is symplectically
perpendicular to $T_x N$.

Let $E$ be the image of the map $(L \g)_\sig^\perp \ra T_xM$ sending
$\eta$ to the fundamental vector field $\eta_M(x)$.  A continuous
inverse is given by the composition of the map $\d_x\Phi:\,E\ra
(L \g)_\sig^0$ with the isomorphism $(L \g)_\sig^0\cong (L \g)_\sig^\perp$. Hence
$E$ is a closed complement to $T_xN$.

By equivariance of the moment map, one has for all $\eta,\zeta\in
L\g$,
$$\omega_x(\eta_M,\zeta_M)=\nu_\xi(\eta_{L\g^*}, 
 \zeta_{L\g^*})=\f{1}{2\pi}\oint\eta \cdot \delta_\xi\zeta. 
$$  
By Lemma \ref{SplitLemma}, $(L \g)_\sig^0$ is a symplectic subspace of
$T_{\Phi(x)}(LG.\Phi(x))$.  It follows that $\d_x\Phi$ restricts to a
symplectic isomorphism from $E$ to $(L \g)_\sig^0$, and also that $T_x N$
is symplectically perpendicular to $E$.

The restriction of $\Phi$ to $Y_\sig$ is a moment map for the
$(LG)_\sig$-action which is equivariant with respect to the affine action
on $(L \g)_\sig^*\oplus\{ 1 \}$ given in Equation
(\ref{AffineAction}). Hence, subtracting $\mu$ gives a moment map
which is equivariant with respect to the usual coadjoint action.
\end{proof}
In particular this shows that Hamiltonian actions of loop groups 
at non-zero level are always proper group actions.

\begin{remark}
If $M$ is a Hamiltonian $L(G\times G)$-manifold with moment map 
$(\Phi_+,\Phi_-)$, we also define cross-sections 
$$ Y_{\sig,-\sig}:=(\Phi_+,\Phi_-)^{-1}(U_\sig\times (-U_\sig)).$$
which are Hamiltonian $(LG)_\sig\times (LG)_\sig$-spaces. 
Note that 
the anti-diagonal action of 
$(LG)_\sig\subset LG$ on $Y_{\sig,-\sig}$ is Hamiltonian with 
moment map the restriction of $\Phi_+ +\Phi_-$, i.e. no shift is 
required. 
\end{remark}

\begin{remark}
In our applications, $M$ will in fact have an invariant  K\"ahler
structure. However, the symplectic cross-sections are usually 
not K\"ahler submanifolds of $M$.
\end{remark}

As in the finite dimensional case, every Hamiltonian $LG$-manifold has
Hamiltonian actions (induced toric flows) of the centers $Z((LG)_\sig)$ on
$LG\cdot Y_\sig$ that commute with the action of $LG$. The moment maps
for these actions are given as the composition of the 
induced toric moment map
$$\ti{\Phi}:=q\circ\Phi:\,M\ra \Alc,$$
with the projection $\t^*\ra \z((L \g)_\sig)^*$.  Here $q:\,L\g^*\ra 
\Alc=L\g^*/LG$ is the quotient map, which can also be written as the
composition of the holonomy map with the quotient map $G\ra G/\Ad(G)$.
As before, if $M = \hat{M} \qu  LG$ is the symplectic reduction of a
Hamiltonian $L(G \times G)$-manifold $\hat{M}$ with moment map
$(\Phi_+,\Phi_-)$ by the anti-diagonal action, the maps
$\ti{\Phi}_\pm$ for $\hat{M}$ descend to a map $\ti{\Phi}:\,M\ra 
\Alc$, whose $\z((L \g)_\sig)^*$-component is a moment map for the induced
action of $Z((LG)_\sig)$ on $\ti{\Phi}^{-1}(U_\sig\cap \Alc)$.  We
refer to these as the residual toric moment map and residual toric flow.

\begin{proposition} \label{DiagonalProp}
Let $H,\,G$ be compact connected simply connected groups and
$\Alc$ the fundamental alcove for $G$.
Let $M$ be a Hamiltonian $L(H\times G \times G)$-Banach manifold,
$M_0 := M \qu LG$ the (possibly singular) reduction by the diagonal
$LG$-action, and $\ti{\Phi}_0: \, M_0 \ra \Alc$ the residual 
toric moment map.  For
every face $\sig$ of $\Alc$ we have a canonical homeomorphism
\begin{equation} \label{FinDiml}
\ti{\Phi}_0^{-1}( U_\sig) \cong Y_{\sig,-\sig} \qu (LG)_\sig.
\end{equation}
If the anti-diagonal action of $LG$-action is free on the zero level
set, then $M_0$ is a smooth Hamiltonian $LH$-Banach manifold and the
above identification is a symplectomorphism.  If $H=\{1\}$ and if
the moment map is proper then $M_0$ is finite dimensional.
\end{proposition}

\begin{proof}
Note that $\ti{\Phi}_0^{-1}( U_\sig)$ is equal to  
the open subset 
$$ LG^2 \cdot Y_{\sig,-\sig} \qu LG \subset M_0. $$
Since $LG^2 \cdot Y_{\sig,-\sig}$ is symplectomorphic to
the symplectic induction
$$\Ind_{(LG^2)_{\sig,\sig}}^{LG^2}(Y_{\sig,-\sig})= LG^2
\times_{(LG^2)_{\sig,\sig}} Y_{\sig,-\sig},$$
the formula (\ref{FinDiml}) follows.  Now suppose that the
anti-diagonal $LG$-action on the zero level set is free. Then the
diagonal action of $(LG)_\sig$ on the zero level set in $Y_{\sig,-\sig}$
is free. Since $(LG)_\sig$ is compact, this implies that $0$ is a regular
value for the moment map for both of these actions, and
(\ref{FinDiml}) is a diffeomorphism of Banach manifolds.  Also, it is
clear that the 2-forms induced by the symplectic form $\omega$ on both
sides of (\ref{FinDiml}) coincide. Since the Meyer-Marsden-Weinstein
theorem holds for Hamiltonian actions of compact groups on Banach
manifolds, it follows that the two-form on $M_0$ is symplectic.
\end{proof}
This shows in particular that the Meyer-Marsden-Weinstein theorem
holds for anti-diagonal $LG$-actions for Hamiltonian 
$L(G\times G)$-actions at
positive level. (It does not hold in general for Hamiltonian actions
of Banach Lie groups.)

\subsubsection{Reductions with respect to coadjoint $LG$-orbits}

Let $G$ be a compact connected simply connected Lie group, with fundamental 
alcove $\Alc$, and let  
$M$ be a Hamiltonian $LG$-Banach manifold. 
For each $\xi\in\Alc$ let
${\O}_{\xi}=LG\cdot\xi$ be the corresponding loop group orbit.
Define the {\bf reduction of $M$ at $\xi$} as the
symplectic reduction by the anti-diagonal $LG$-action
$$ M_{\xi}=M\times\O_{*\xi}\qu LG. $$
Letting $\sigma$ be a face of $\Alc$ such that
$\xi\in U_{\sig}$ and $Y_{\sigma}$ the
corresponding symplectic cross-section of $M$, one also has
\begin{equation} 
M_{\xi}=Y_{\sigma}\times 
{\O}'_{*\xi}\qu (LG)_\sig,
\label{FiniteReduction}
\end{equation}
where $\O_{*\xi}' = \O_{*\xi} \cap ((L \g)_\sig^*\times\{1\})=
\widehat{(LG)_\sig}\cdot *\xi$ is the orbit for the compact group
$\widehat{(LG)_\sig}$.  From (\ref{FiniteReduction}), one finds in
particular that if the moment map is proper then
$M_{\xi}$ is a finite dimensional symplectic quotient. 
In the case that the quotient is singular, the local structure of 
its singularities can be described by normal form theorems as 
in Sjamaar-Lerman \cite{SL}. 

If $M$ is a Hamiltonian $LG$-manifold at level $\levi\in \N$, $L \ra M$ an
$\widehat{LG}$-equivariant pre-quantum line bundle, and
$\xi \in \levi \Alc \cap \Lambda^*$ then we define
\begin{equation} \label{LineReduction}
L_{\xi} := L \boxtimes \Xi(\O_{*\xi})\qu LG, 
\end{equation}
where $\Xi(\O_{*\xi})$ is the pre-quantum line bundle of $\O_{*\xi}$
as in example \ref{CoadjointOrbits}.
If the anti-diagonal action of $LG$ is locally free on the zero level set then
$L_{\xi}$ is a pre-quantum (orbifold) line bundle over
$M_{\xi}$.

\subsection{Convexity theorems for Hamiltonian loop group actions} 
In the following theorem, we use symplectic cross-sections 
to derive convexity and connectedness properties for the 
moment map of Hamiltonian $LG$-manifolds.  Some of these results
have been proved by Chang \cite{C2} using similar methods.

\begin{theorem} \label{FiniteDimensional}
Let $G$ be a compact connected simply connected Lie group, and 
$M$ a connected Hamiltonian $LG$-manifold with proper moment map
$\Phi:M \ra L \g^*$.
\begin{enumerate}
\item
For any face $\sig$ of the fundamental alcove $\Alc$ such that
$\Phi(M)\cap\sig\not=\emptyset$, the corresponding symplectic
cross-section $Y_\sig$ is finite dimensional and connected.
\item The fibers of $\Phi$ are connected. 
\item The intersection $\Phi(M)\cap\Alc$ is 
      a convex polytope. 
\end{enumerate}
\end{theorem}

\begin{proof}
Finite dimensionality follows from the fact that the restriction of 
$\Phi$ to the submanifold $Y_\sig=\Phi^{-1}( U_\sig)$ is proper 
as a map to $U_\sig$. Since $U_\sig$ is finite dimensional, this 
is only possible if $Y_\sig$ is finite dimensional. 

We next prove that all non-empty cross-sections $Y_\sig$ are
connected. Since the coadjoint orbits $LG/(LG)_\sig$ are simply
connected, this is the case if and only if all flow-outs $LG\cdot
Y_\sig=LG\times_{(LG)_\sig} Y_\sig$ are connected.  By properness of
$\Phi$, the number of connected components $Y_\sig^i$ of $Y_\sig$ is
finite. The fact that $LG \cdot Y_\sig$ is connected follows once 
we can show  the transitivity property 
\begin{equation} \label{Transitivity}
LG\cdot(Y_\sig^i\cap Y_\tau^j)\not=\emptyset,\ \ 
LG\cdot(Y_\tau^j\cap Y_\kappa^k)\not=\emptyset\ 
\Longrightarrow 
LG\cdot(Y_\sig^i\cap Y_\kappa^k)\not=\emptyset, 
\end{equation}
since the collection of all $LG\cdot Y_\sig^i$ is a finite open
covering of $M$, and any two points in $M$ can be joined by a path.

To show (\ref{Transitivity}), let $x\in Y_\sig^i\cap Y_\tau^j$ and
$y\in Y_\tau^j\cap Y_\kappa^k$.  The fact that the restriction of
$\ti{\Phi}$ to $Y_\tau^j$ is proper as a map into $U_\tau$ implies
that it has connected fibers and that the image
$\ti{\Phi}(Y_\tau^j)\subset \Alc$ is convex.  This follows from the
Condevaux-Delzord-Molino technique as explained in \cite[p. 29]{FR}
and \cite{HNP}, or by using symplectic cutting \cite[Remark
5.2]{LMTW}.  In particular, the line segment from
$\alpha=\ti{\Phi}(x)$ to $\beta=\ti{\Phi}(y)$ is contained in
$\ti{\Phi}(Y_\tau^j)$.  It follows that there exists a path $\gamma$
in $LG\cdot Y_\tau^j$ connecting $x$ and $y$ whose image under
$\ti{\Phi}$ is the line segment $\ol{\alpha\beta}$.  The interior of
$\ol{\alpha\beta}$ intersects $U_\sig \cap \Alc$ and
$U_\kappa \cap \Alc$, and is therefore contained in $U_\sig \cap
U_\kappa \cap \Alc$.  It follows that $\gamma^{-1}
(\text{int}(\ol{\alpha\beta}))\subset LG\cdot(Y_\sig^i\cap
Y_\kappa^k)$.

Since all cross-sections $Y_\sig$ are connected and all restrictions
$\Phi|Y_\sig$ have connected fibers \cite{LMTW}, it follows that
$\Phi$ has connected fibers.  Since $\Phi$ is proper, the number of
$(LG)_\sig$-conjugacy classes of stabilizer groups (orbit types) for the
$(LG)_\sig$-action on $Y_\sig$ is finite. By \cite[Remark 5.2]{LMTW} or
the Condevaux-Delzord-Molino technique, this implies that
$\ti{\Phi}(Y_\sig)$ is the intersection of $U_\sig\cap\Alc$ with a
convex polyhedron. Since $M$ is connected, and since a closed
connected set is convex if and only if it is locally convex, this
shows that the image $\ti{\Phi}(M)$ is a convex polytope.
\end{proof}

\begin{remark}  
As far as we know the Atiyah-Pressley convexity Theorem \cite{AP},
which is a Kostant-type convexity theorem on projections of coadjoint
loop group orbits, does not fit into this framework. Note that
the relevant moment map, given by projection of the orbit to
$\t^*$ together with the energy function, is not proper.
\end{remark}

\begin{corollary}\label{ToricConvexity}
(Convexity properties of the induced toric moment map.)  
Let $H,\,G$ be compact connected simply connected Lie groups, with fundamental 
alcoves $\lie{B},\,\Alc$, and $M$ a Hamiltonian 
$L(H\times G\times G)$-manifold with proper moment map. 
Let $M_0$ be the symplectic
reduction with respect to the anti-diagonal $LG$-action, and
$\ti{\Phi}_0:\,M_0\ra \lie{B}\times \Alc$ the product of
the induced toric moment map for the induced $LH$-action and the residual 
toric moment map
obtained from the $ \{e\}\times\{e\} \times LG$-action on $M$. Then
the image of $\ti{\Phi}_0$ is a convex polytope.
\end{corollary}

\begin{proof} 
Let $\Phi:\,M\ra L(\h\oplus \g\oplus \g)^*$ be the moment map
for the $L(H\times G\times G)$-action, and $\ti{\Phi}:\,M\ra 
\lie{B}\times \Alc\times \Alc$ the corresponding
toric moment map. Consider the involution $\kappa$ of $\,
\lie{B}\times \Alc\times \Alc$ defined by
$$\kappa(x_1,x_2,x_3)=(x_1,* x_3,* x_2)$$
where $*$ is the Cartan involution given by $*\mu = q(-\mu)$.  The
image of the induced map $\ti{\Phi}_0$ is given by the image of the
intersection $\ti{\Phi}(M)\cap ( \lie{B}\times \Alc\times
\Alc)^\kappa$, where $( \lie{B}\times \Alc\times \Alc)^\kappa$
denotes the fixed point set of the involution $\kappa$, under the
projection of $ \lie{B}\times \Alc\times \Alc$ to the first and
second factor.
\end{proof}

\subsection{Application to Yang-Mills theory over a Riemann surface}

Let $\Sig$ be a compact, connected, oriented Riemann surface of genus
$g$ with $b$ boundary components, and $G$ a compact, connected and 
simply connected Lie group, with fundamental alcove $\Alc$.
  
\subsubsection{Holonomy description}
As in the case $\p\Sig=\emptyset$, the moduli space $\M(\Sig)$ 
of flat $G$-connections on $\Sig\times G$
admits an alternative description in terms of holonomies.   
\begin{theorem} \label{HolonomyDescription}
If $b\ge 1$, 
the  moduli space $\M(\Sigma)$ is isomorphic to the set of 
$ (a,c,\xi) \in 
G^{2g}\times G^{b-1}\times (L\g^*)^b $ such that 
$$
\prod_{i=1}^{2g}[a_{2i-1},a_{2i}]= \prod_{i=1}^b \, \Ad_{c_i} \, \Hol(\xi_i)
$$
where $c_1 = 1$. This is a smooth submanifold of 
$G^{2g}\times G^{b-1}\times (L\g^*)^b$ 
and the identification with  $\M(\Sigma)$
is an $LG^b$-equivariant diffeomorphism. Here the action of 
$g=(g_1,\ldots, g_b)\in LG^b$ 
on $G^{2g}\times G^{b-1}\times (L\g^*)^b $ is given by
$$ g\cdot a_i=\Ad_{g_1(0)}a_i,\ \ g\cdot c_j=g_1(0)\,c_j g_j(0)^{-1}
,\ \ g\cdot \xi_j=\Ad_{g_j}\cdot\xi_j-(\,g_j)'\,{g_j}^{-1}
$$
and the moment map is given by projection to the $(L\g^*)^b$-factor.
\end{theorem} 
\begin{proof}
The diffeomorphism is given as follows.  Let $B_1,\ldots,B_b$ be the
boundary circles of $\Sig$ equipped with fixed parametrizations $B_j
\cong S^1$ compatible with the orientations on $B_j$.  This gives
identifications $\G(\p \Sig) \cong LG^b$ and $\Omega^1(\p \Sig,\g)
\cong (L \g^*)^b$ and base points $x_j \in B_j$.  Now let
$\rho_1,\ldots,\rho_{2g}$ be smooth loops based at $x_1$, and
$\gamma_2,\ldots,\gamma_b$ smooth paths from $x_1$ to
$x_2,\ldots,x_n$, in such a way that $\pi_1(\Sig)$ is generated by the
$\rho_i$ together with $B_1$ and $\gamma_j^{-1}\, B_j\, \gamma_j$,
subject to the relation
$$ \prod_{i=1}^g [\rho_{2i-1},\rho_{2i}]=B_1\, \big(\gamma_2^{-1}\,
B_2\, \gamma_2\big) \ldots \big(\gamma_b^{-1}\, B_b\, \gamma_b\big).$$
Consider the map 
$$\hat{f}:\,\A_F(\Sig)\ra G^{2g}\times G^{b-1}\times (L\g^*)^b$$
that takes any flat connection $A$ to $(a,c,\xi)$
where $a_i$ is the 
holonomy around $\rho_i$, $c_j$ the parallel transport along 
$\gamma_j$, and $\xi_j$ the restriction of $A$ to the boundary loop $B_j$.
Then 
\begin{equation}
 \prod_{i=1}^g [a_{2i-1},a_{2i}] = \prod_{j=1}^b \Ad_{c_j}(\Hol(\xi_j)),
\label{HolonomyCondition}
\end{equation}
where we set $c_1:=1$.
The map
$\hat{f}$ is equivariant with respect to the action of the gauge group 
$\G(\Sig)$ given on $G^{2g}\times G^{b-1}\times (L\g^*)^b$ by 
$$ g\cdot a_i=\Ad_{g(x_1)}a_i,\ \ g\cdot c_j=g(x_1)\,c_j g(x_j)^{-1}
,\ \ g\cdot \xi_j=(g|B_j)\cdot\xi_j. $$
Moreover, $\hat{f}$ is surjective onto the set of all $(a,c,\xi)$
satisfying (\ref{HolonomyCondition}): First, as in the case without
boundary, one can construct a smooth connection that has the required
holonomies $a_i$ and $\Ad_{c_j}(\Hol(\xi_j))$.  Secondly, one can act
by an element $g$ of $\G(\Sig)$ with $g(x_1)=1$ (which does not change
the holonomies) to obtain the required values of $\xi_j$.  This gives
the correct values for the parallel transport $c_j$ along the curves
$\gamma_j$, up to an action of the centralizer of $\Hol(\xi_j)$ from
the right. Now evaluation at $x_j$ gives an isomorphism of
$Z(\Hol(\xi_j))$ with the stabilizer $(LG)_{\xi_j}$.  Thus finally we
can act by an element $g'\in\G(\Sig)$ with $g'(x_1)=1$ and $g'(x_j)\in
Z(\Hol(\xi_j))$ for $j\ge 2$ to obtain the correct values of $c_j$.

We next show that the fiber of $\hat{f}$ over $(a,c,\xi)$ is 
equal to $\G_\p(\Sig)$. To see this, one may assume after acting by 
$\G(\Sig)$ that $\xi$ is smooth. Let $A_1,\,A_2\in \hat{f}^{-1}(a,c,\xi)$. 
For $x\in \Sig$, choose a smooth path $\lambda$ from $x_1$ to $x$, 
and let $g(x)\in G$ defined by parallel transport along $\lambda$ by 
$A_1$, followed by parallel transport along $\lambda^{-1}$ by $A_2$. 
This is independent of the choice of $\lambda$, and therefore gives 
a smooth function $g:\,\Sigma\ra G$ with $g\in\G_\p(\Sig)$ 
and $g\cdot A_1=A_2$. By a similar argument, one shows that 
the kernel of the tangent map to $\hat{f}$ is equal to the 
tangent space to the $\G_\p(\Sig)$-orbit. 
Hence $\hat{f}$ descends to a smooth embedding
$$ f:\,\M(\Sig) \ra G^{2g}\times G^{b-1}\times (L\g^*)^b$$
with image equal to the set (\ref{HolonomyCondition}). 
\end{proof}

It is immediate from this description that the moment map is proper.
Equivalently, the cross-sections are finite-dimensional.  This fact
can also be proved in the gauge theory description by writing the
tangent space to the cross-section as the solution space to an
elliptic boundary value problem.  Since $G$ is a product of simple
groups, thus $(G,G)=G$, the holonomy description also shows that for
$g\ge 1$ the moment map is surjective onto $(L\g^*)^b$ so that the
moment polytope is simply $\Alc^b$ and the convexity Theorem
\ref{FiniteDimensional} is vacuous.  On the other hand, in the case
$g=0$ of a $b$-punctured sphere, (\ref{HolonomyCondition}) shows that
the moment polytope is equal to the image of the subset
$$ \big\{d\in G^b\big|\,\prod_{j=1}^b d_j\,=1\big\}$$
under the quotient map $G^b\ra \Alc^b$. 

\begin{example}\label{ModuliTwoPuncturedSphere2}
By (\ref{HolonomyCondition}), the moduli space $\M(\Sig^2_0)$ of the
two-punctured sphere (annulus) is equivariantly diffeomorphic to the
subset $G\times L\g^*\times L\g^*\ni(c,\xi_1,\xi_2)$ defined by
$\Hol(\xi_1)\,\Ad_c\,\Hol(\xi_2)=1$. Consider the Hamiltonian
$L(G\times G)$-manifold $X\cong LG\times L\g^*$ introduced in Example
\ref{ModuliTwoPuncturedSphere}. The map
$$X\ra \M(\Sig^2_0),\,(g,\eta)\mapsto (g(0),\,g\cdot
\eta,\,-\eta) $$
is an equivariant diffeomorphism, preserving the 
moment maps. Since $\M(\Sig^2_0)$ and $X$ are multiplicity free, i.e. 
since all reduced spaces $X_{\xi_1,\xi_2}$ are points, this map 
is necessarily a symplectomorphism.  
\end{example}

Along similar lines, one can prove the holonomy description for 
the moduli spaces $\M(\Sig)_{\xi_1,\ldots,\xi_b}$ (see e.g. 
\cite{AM}): 
\begin{theorem}
Let $\xi_1,\ldots,\xi_b \in \Alc$, and
$\mathcal{C}_1,\,\ldots,\,\mathcal{C}_b\subset G$ the corresponding conjugacy
classes.  The moduli space $\M(\Sig)_{\xi_1,\ldots,\xi_b}$ is
homeomorphic to the the quotient of the subset $ (a,d)\in
G^{2g}\times\mathcal{C}_1\times\,\ldots\times\,\mathcal{C}_b $ given
by
\begin{equation}
 \prod_{i=1}^g [a_{2i-1},a_{2i}] = \prod_{j=1}^b d_j,
\end{equation}
by the conjugation action of $G$.
\end{theorem}

\subsubsection{Symplectic cross-sections}

Consider once again the Hamiltonian $LG^b$-manifold $\M(\Sig)$ for $b
> 0$. In the holonomy description,  
the cross-sections $Y_{\sig_1\ldots\sig_b}$ for faces 
$\sig=(\sig_1,\ldots,\sig_b)$ of $\Alc^b$ 
are the smooth
submanifolds of $G^{2g}\times G^{b-1}\times (L\g^*)^{b}$ given by the
condition (\ref{HolonomyCondition}) together with the requirement
$\xi_j\in U_{\sig_j}$. For the dimension of $Y_{\sig_1\ldots\sig_j}$,
we find
\begin{eqnarray*}
\dim(Y_{\sig_1\ldots\sig_j})&=&(2g+b-1)\dim G+\sum \dim U_{\sig_j}
-\dim G\\
&=&(2g-2)\dim G  +\sum(\dim G+\dim (LG)_{\sig_j}).
\end{eqnarray*}
The extended moduli spaces of Chang \cite{C1}, Huebschmann \cite{H}
and Jeffrey \cite{J} are defined, for a Riemann surface with a single
boundary component, as the subset of $G^{2g}\times\g^*$ given by the
condition $\prod [a_{2i-1},a_{2i}]=\exp(2\pi\xi)$.  On a neighborhood
of $\xi=0$, this is a smooth submanifold and can be given a symplectic
structure; however for larger $\xi$ one encounters singularities and
degeneracies of the symplectic form.  One can view the extended moduli
space as the symplectic cross-section corresponding to $\sigma=\{0\}$.
The degeneracies find a natural explanation in the full space $\g^*$
not being a slice for $LG$, whence the cross-section may not be chosen
too big. In \cite{J}, more general extended moduli spaces are defined
for all central elements $c\in Z(G)$ by the modified condition $\prod
[a_{2i-1},a_{2i}]=c\,\exp(2\pi\xi)$; these may be identified, using
Remark \ref{CentralElements}, with the cross-sections corresponding to
$\sig=\{\eta\}$, where $\eta\in\Alc$ is the unique vertex such that
$\exp(-2\pi\eta)=c$.

\subsubsection{Pants decomposition and Goldman twist flows}
\label{GoldmanSec}

For every compact connected Riemann surface $\Sig$ of genus $g$ with
$b$ boundary components, except for the cases
$(g,b)=(0,0),\,(0,1),\,(0,2),\,(1,0)$, there exist embedded circles
$C_1,\ldots,C_r\subset \Sig$ that decompose $\Sig$ into a union of
3-holed spheres (pairs of pants) $\Sig_1,\ldots \Sig_l$. For the
number $l$ of pants and $r$ of separating curves in this
decomposition, one has $3l=b+3r$, where:
 \begin{eqnarray*}
r=3g+b-6 &\text{if}  & b\ge 1,\,g\ge 2,\\
r=3g-3 &\text{if} & b=0,\,g\ge 2,\\
r=b  &\text{if}  & b\ge 1,\,g=1,\\
r=b-1 &\text{if}  & b\ge 3,\,g=0.
 \end{eqnarray*}
By repeated application of Theorem \ref{GluingEqualsReduction}, the
moduli space $\M(\Sig)$ is obtained as a symplectic reduction
$$  \M(\Sig)=\M(\Sig_1)\times\ldots\times \M(\Sig_l)\qu LG^r.$$
As in the finite dimensional case, this induces a residual toric 
flow on $\M(\Sig)$, i.e. a densely defined, effective Hamiltonian 
$T^{r+b}$-action which commutes with the action of $\G(\p\Sig)$, 
known as the {\bf Goldman twist flow}.
The corresponding residual toric moment map, called here the Goldman map, is
$$ {\lie{G}}: \, \M(\Sig) \ra \Alc^{r+b}, \, [A] \ra ( G \cdot \Hol_{B_1}(A),
\ldots , G \cdot \Hol_{B_b}(A),\,  G \cdot \Hol_{C_1}(A),
\ldots , G \cdot \Hol_{C_r}(A) ).$$
As a special case of  Corollary \ref{ToricConvexity}, we have:
\begin{theorem}
The image ${\lie{G}}(\M(\Sig))\subset \Alc^{b+r}$ of the Goldman map
${\lie{G}}$ is a $(b+r) \cdot \dim T$-dimensional convex polytope.
\end{theorem}

This polytope is given explicitly as follows.  Let $P \subset \Alc^3$
denote the polytope which is the moment polytope for
the three-holed sphere.  For each $\Sig_j$ let $B_j^1,B_j^2,B_j^3
\subset \p \Sig_j$ denote its boundary components, and let 
$\Alc_j^\nu$ be a copy of the fundamental alcove associated to 
$B_j^\nu$. 
Define an involution
$ \kappa: \Alc^{3l} \ra \Alc^{3l} $ by
$$( \kappa(x) )^\nu_j =\left\{
\begin{array}{lc}
 x^\mu_l  &  \hbox{ \ if \ } B^\nu_j = B^\mu_l \hbox{\ for some $(l,\mu)$ } \\
 x^\nu_j  &  \hbox{ \ otherwise. } 
\end{array}\right.
$$
Then $\lie{G}(\M(\Sig))$ is the image of the intersection $P^l
\cap (\Alc^{3l})^\kappa$ where $ (\Alc^{3l})^\kappa$ denotes the fixed
point set of the involution $\kappa$, under the projection $\Alc^{3l}
\ra \Alc^{b + r}$.

\begin{remark}
For the case $G = SU(2)$ and $\partial \Sig =\emptyset$ the dimension
$(2g-2)\dim G$ of the moduli space $\M(\Sig)$ is precisely twice the 
dimension of this polytope, which means that the Goldman flow  
gives a completely integrable torus action
on a dense subset of $\M(\Sig)$. This has been studied extensively 
in Jeffrey-Weitsman \cite{JW}.
\end{remark}

By Proposition \ref{FiniteDimensional} if $\sig$ is a face of the
fundamental alcove $\Alc$ then we can write $\lie{G}^{-1}(U_\sig
\cap \Alc) \subset \M(\Sig)$ as a symplectic reduction of
finite-dimensional symplectic manifolds:

\begin{lemma} 
Let $\hat{\Sig}$ be a compact oriented Riemann surface (possibly
disconnected) and $\Sig$ the Riemann surface formed by gluing along
two boundary components $B_\pm \subset \partial \hat{\Sig}$,
and let $\Phi_\pm : \M(\hat{\Sig}) \ra L \g^*$ denote the moment
maps for the $LG$-actions corresponding to $B_\pm$.  
Then there is a symplectic diffeomorphism
$$\lie{G}^{-1}(U_\sig \cap \Alc) = \Phinv_+(U_\sig) \cap
\Phinv_-( - U_\sig) \qu  (LG)_\sig .$$
\end{lemma}

\section{Quantization}

In this section, we define the quantization of moduli spaces of flat
connections using the equivariant index (Riemann-Roch number) of
$\text{Spin}_c$-Dirac operators, and introduce techniques for dealing
with singular quotients.

\subsection{Riemann-Roch numbers for almost complex orbifolds}
\label{RiemannRochNumbers}
Let $M$ be a compact almost complex $G$-manifold. Then $M$ has a
canonical $G$-invariant $\text{Spin}_c$-structure. Given a Hermitian
orbifold vector-bundle $E\ra M$, one can consider the G-equivariant
$\Spin_c$-Dirac operator $\dirac_E$ corresponding to this
$\Spin_c$-structure, twisted by $E$ (see e.g. \cite{LM}).  The
equivariant index of $\dirac_E$ is a virtual representation of $G$
which we call the {\bf equivariant Riemann-Roch number},
$$\chi=\RR(M,E):=\text{ind}_G(\dirac_E)\in \Rep(G),$$
where $\Rep(G)$ denotes the representation ring for $G$.  An explicit
expression for $\RR(M,E)$ is given by the Equivariant Index Theorem of
Atiyah-Segal-Singer.

If $M$ is an almost complex $G$-orbifold, and $E\ra M$ a $G$-orbifold
vector bundle\footnote{An orbifold vector bundle is modeled locally on
the quotient of a smooth vector bundle by a finite group of bundle
automorphisms.  We will usually refer to these simply as vector
bundles.} the same definition applies. A formula for the character is
given by the orbifold index theorem of Kawasaki \cite{Ka} (the
equivariant version of which was proved by Vergne \cite{Vg2}).
Riemann-Roch numbers satisfy the following functorial properties:
\begin{enumerate}
\item (Products)
      If $E_i\ra M_i$ are $G_i$-equivariant vector bundles over compact
      almost complex $G_i$-orbifolds ($i=1,2$), then 
       $$\RR(M_1\times M_2,E_1\boxtimes E_2)=\RR(M_1,E_1)\otimes
       \RR(M_2,E_2).$$ 
\item (Conjugation)
      Let $E\ra M$ be a $G$-equivariant orbifold bundle over a 
      compact almost complex $G$-orbifold. Let $M^*$ denote $M$ with 
      the opposite almost complex structure. 
      Then 
      $$ \RR(M^*,E^*)=\RR(M,E)^*, $$
      where, for all $\chi\in \Rep(G)$, $\chi^*\in \Rep(G)$ denotes the dual 
      representation.
\item (Induction) Let $\sigma\subset\t^*_+$ be an open face of the 
       positive Weyl chamber, with corresponding stabilizer group $G_\sig$.
       Let $\Lambda^*_{\sig,+}\supset \Lambda^*_+$ be the dominant weights 
       for $G_\sig$, and 
       $$\Ind_{G_\sig}^G:\,\Rep(G_\sig)\ra \Rep(G)$$ 
       denote
       holomorphic induction.  Let $G/G_\sig$ be equipped with the
       standard complex structure coming from its identification with
       a coadjoint orbit $G\cdot\mu$, for any $\mu\in\sig$.  For every
       almost complex $G_\sig$-orbifold $Y_\sig$, the associated
       bundle $G\times_{G_\sig}Y_\sig$ has a naturally induced almost
       complex structure. Given a $G_\sig$-invariant orbifold vector
       bundle $E_\sig\ra Y_\sig$, we have
$$ \RR(G\times_{G_\sig}Y_\sig,G\times_{G_\sig}E_\sig)
=\Ind_{G_\sig}^G \RR(Y_\sig,E_\sig). $$
For a proof, see e.g. \cite{M2}.
\end{enumerate}

If $M$ is a Hamiltonian $G$-orbifold with moment map $\Phi:\, M \ra
\g^*$, and $E\ra M$ a $G$-equivariant vector bundle, one can always
choose a $G$-invariant compatible almost complex structure $J$ on $M$
to define $\RR(M,E)$. Since any two $J$'s are equivariantly homotopic,
this definition does not depend on the choice of $J$.  We will be
mostly interested in the case that $E = L$ is a {\bf pre-quantum line
bundle}, i.e., $L$ is a $G$-equivariant Hermitian line bundle equipped
with invariant connection whose curvature is $- 2\pi i$ times the
symplectic form, and the fundamental vector fields on $L$ and $M$ are
related by the Kostant formula
\begin{equation} 
\xi_L=\text{Lift}(\xi_M)+2\pi\l\Phi,\xi\r \,\f{\p}{\p\phi} 
\label{LiftFormula}
\end{equation}
where $\f{\p}{\p\phi}$ is the generating vector field for the scalar
circle action on the fiber.

\subsection{Desingularization for quotients by compact groups}
\label{SymplecticQuotients}
In this section we define Riemann-Roch numbers for singular symplectic
quotients by compact groups and state the ``quantization commutes 
with reduction''
Theorem, referring to  \cite{MS} for more details. 

Suppose that $M$ is a Hamiltonian $H\times G$-orbifold
($H,G$ compact), such that the $G$-action has proper moment map
$\Phi:\,M\ra \g^*$, and let $L\ra M$ be a $H \times G$-equivariant
pre-quantum bundle. Let $M_0=M\qu G$ be the symplectic quotient by $G$
and $ L\qu G = L_0:=L|\Phi^{-1}(0)/G $ the quotient bundle.

It can be shown (\cite{MS}) that if $\Phi^{-1}(0)$ is contained
in a fixed infinitesimal $G$-orbit type stratum (that is, if the
dimension of the isotropy group $G_x$ does not jump as $x$ varies in
$\Phi^{-1}(0)$) the level set $\Phi^{-1}(0)$ is an orbifold
and consequently $M_0$ is a symplectic orbifold and $L_0$ with its 
induced connection an orbifold pre-quantum
line bundle.

If $\Phi^{-1}(0)$ and $M_0$ have more serious singularities, we define
$\RR(M_0,L_0)$ as the $H$-equivariant Riemann-Roch number of an
orbifold line bundle $\ti{L}_0\ra \ti{M}_0$ obtained from $L_0\ra
M_0$ by means of {\bf partial desingularization}.  One procedure to
obtain a desingularization is due to Kirwan \cite{Ki}.  Roughly
speaking, this is an inductive procedure involving a sequence of
symplectic blow-ups on $M$, resulting in a new Hamiltonian $H \times
G$-orbifold $\ti{M}$, followed by a symplectic quotient
$\ti{M}_0=\ti{M}\qu G$.  Since blow-ups in the symplectic category
depend on both the choice of Darboux charts around the blow-up locus
and on the ``size'' of the blow-up corresponding to the cohomology
class of the symplectic form on the exceptional divisor (see
\cite{McS} for details), the partial desingularization $\ti{M}_0$
defined this way is unique only up to {\bf symplectic homotopy}, and
the bundle $\ti{L}_0\ra \ti{M}_0$ is unique up to isomorphism. This,
however, is enough to obtain a well-defined Riemann-Roch number
$\RR(\ti{M}_0,\ti{L}_0)$.

A simpler, but less canonical way of desingularizing a symplectic quotient 
is to shift the value of the moment slightly.  
For all $\alpha\in\g^*$, let
$M_\alpha=\Phi^{-1}(\alpha)/G_\alpha$
(with the induced $H$-action) be the symplectic quotient.  For
$\alpha\in\Phi(M)$ close to zero and generic\footnote{We call a value
$\alpha$ of the moment map $\Phi$ {\bf generic} if the restriction of
$\Phi$ to $\Phi^{-1}(\alpha)$ has maximal rank. This implies that
$\Phi^{-1}(\alpha)$ is a smooth suborbifold and that the dimension of
$G_x$ for $x\in \Phi^{-1}(\alpha)$ is constant.}, one can take $ (L\,
| \, \Phi^{-1}(\alpha))/G_\alpha \ra M_\alpha$ as a desingularization
of $L_0\ra M_0$.

One should keep in mind that the diffeotype of the desingularization
obtained in this way will in general depend on the choice of the
shift.  In fact, if $G$ non-abelian, the two spaces may even have
different dimension (so it is perhaps misleading to call $M_{\alpha}$
a desingularization of $M_0$ in this case.)  It turns out that the
desingularization procedures described above give the same
Riemann-Roch numbers:

\begin{proposition}\cite{MS} Suppose $(M,\omega)$ is 
a Hamiltonian $H\times G$-orbifold, such that the $G$-action has
proper moment map $\Phi:\,M\ra\g^*$, and let $L\ra M$ be an $H\times
G$-equivariant pre-quantum line bundle. Then there is a neighborhood $U$
of $0\in\g^*$ such that for all generic $\alpha\in U\cap\Phi(M)$, one
has $\RR(M_0,L_0)=RR(M_\alpha, (L\, | \, \Phi^{-1}(\alpha))/G_\alpha)$
as representations of $H$, where the first $\RR$-number is defined by
Kirwan's partial desingularization.
\end{proposition}

\begin{theorem}[Quantization commutes with reduction] 
\label{QuantizationReduction} \cite{M2,MS} 
Let $(M,\om)$ be a compact Hamiltonian $G \times H$-orbifold and $L
\ra M$ a $G \times H$-equivariant pre-quantum line bundle.  Then the
$G$-invariant part of $RR(M,L)$ equals the Riemann-Roch number of the
symplectic quotient:
$$ RR(M,L)^G = RR(M \qu  G,L \qu  G) $$
where the right-hand side is defined via partial desingularization
if $0$ is not a regular value of the $G$-moment map.
\end{theorem}

If $\mu\in\Lambda^*_+$ is a dominant weight, the coadjoint orbit 
$\O_\mu=G\cdot\mu$ with its KKS form has a pre-quantum line 
bundle $\Xi(\O_\mu)\ra \O_\mu$.  Given a $G$-equivariant pre-quantum line
bundle $L\ra M$ as above, we define the reduced pre-quantum line bundle by 
$$ L_\mu := (L\boxtimes \Xi(\O_{*\mu}))\qu G\ra M_\mu=(M\times 
{\O}_{*\mu})\qu G.$$
As a direct corollary to Theorem \ref{QuantizationReduction}, the 
multiplicity $N(\mu)$ for the weight  $\mu$ to occur in $\RR(M,L)$ 
is given by $\RR(M_\mu,L_\mu)$, defined by canonical desingularization 
if necessary.

\subsection{Desingularizations for loop group quotients}
\label{DesingSec}

We now define desingularizations for reductions by loop group actions.
Let $G$ be a compact connected simply-connected Lie group and $M$ a
Hamiltonian $L(G\times G)$-manifold with proper moment map at level
$\levi$ and $L(G\times G)$-equivariant pre-quantum line bundle $L$.  We
will show that the reduction $L \qu LG\to M\qu LG$ by the anti-diagonal action
can be written in a canonical way as a reduction in finite dimensions,
so the desingularization of $L \qu LG$ can be carried out as in the
previous section.

Let $\M^\levi(\Sig^3_0)$ denote the moduli space for the three-holed sphere 
at level $\levi$,
$L^\levi(\Sig^3_0)$ its pre-quantum line bundle and 
$$M^{(1)} = M \times \M^\levi(\Sig^3_0) \qu LG\times LG,\ \ \ \
L^{(1)} = L \times L^\levi(\Sig^3_0) \qu LG\times LG$$
the quotient by the product of the two anti-diagonal $LG$-actions,
defined by pairing each $LG$-factor for $M$ with an $LG$-factor for
$\M^\levi(\Sig^3_0)$.  Since these actions are free, $M^{(1)}$ is a
Hamiltonian $LG$-manifold and $L^{(1)}$ an $\widehat{LG}$-equivariant
pre-quantum line bundle.  Because the reduction of
$\M^\levi(\Sig^3_0)$ at $0$ is the moduli space $\M^\levi(\Sig^2_0)$,
the reduction of $M^{(1)}$ at $0$ is given by 
$$M^{(1)}_0 = (M \times \M^\levi(\Sig^2_0) \qu LG) \qu LG = M \qu LG $$
by Example \ref{ModuliTwoPuncturedSphere} and
\ref{ModuliTwoPuncturedSphere2}.  Let $Y^{(1)}_{\{ 0 \}}$ be the
symplectic cross-section for $M^{(1)}$ at $\{ 0 \} \subset \Alc$.
Then 
\begin{equation} \label{AddPuncture}
 M \qu LG \cong M^{(1)}_0 \cong Y^{(1)}_{\{ 0 \}} \qu G, \ \ \ \
   L\qu LG\cong L^{(1)}_0\cong L^{(1)}_0|\,Y^{(1)}_{\{ 0 \}}\qu G
\end{equation}
which proves the claim.   
In particular, the Riemann-Roch numbers of the line bundles
in (\ref{LineReduction}) can be defined by desingularization.
Note that in the case of moduli spaces, the above procedure
corresponds to ``adding a puncture'' to the Riemann surface.

\section{Symplectic Surgery}

Recall from Subsection \ref{GoldmanSec} that if $\Sig$ is a Riemann
surface with $b$ boundary components and $\xi_1,\ldots,\xi_b \in
\Alc$, the Goldman functions corresponding to a pants decomposition of
$\Sig$ generate a densely defined torus action on the moduli space
$\M(\Sig)_{\xi_1,\ldots,\xi_b}$.  In this section, we describe how to
use this torus action in conjunction with Lerman's symplectic cutting
procedure to decompose the moduli space into pieces which can be
expressed as symplectic reductions in finite dimensions.

\subsection{Symplectic cutting}\label{SymplecticCutting}
We briefly recall Lerman's cutting construction \cite{L}. Let 
$M$ be a Hamiltonian $S^1$-orbifold, with moment map $\psi:\,M\ra \R$.
Consider the action of $S^1$ on the product $M\times\C^-$
given by $e^{i\phi}\cdot(m,z)=( e^{i\phi}\cdot m ,\,e^{-i\phi}\, z)$,
with moment map 
$$\ti{\psi}(m,z)=\psi(m)-|z|^2 .$$
The zero level set of $\ti{\psi}$ is a union of $\psi^{-1}(0)\times
\{0\}$ and the set of all $(m,z)$ with $\psi(m)=|z|^2>0$. 
Suppose that $0$ is a regular value of $\psi$. Then $0$ is also a regular 
value of $\ti{\psi}$, and the reduced space $M_+:=M\times\C^-\qu S^1$ 
is a symplectic orbifold. As a topological space, $M_+$ is
obtained from the manifold with boundary $\psi^{-1}(\R_{\ge 0})$
by collapsing the boundary by 
the nullfoliation of the pullback of the symplectic form. 

It was shown by Lerman that this is also true
symplectically:

\begin{proposition}[Lerman] 
Let
$$M_+:=(M\times\C)\qu S^1$$
be the {\bf cut space}. The canonical maps
$$
\iota_0:\,M_0\ra M_+,\ \ \iota_{>0}:\,\psi^{-1}(\R_{>0})\ra M_+ 
$$
are smooth symplectic embeddings, and the normal bundle of $M_0$ in $M_+$ is
canonically isomorphic to the associated bundle
$\psi^{-1}(0)\times_{S^1}\C$, where $S^1$ acts on $\C$ with weight $-1$. 
Given a Hamiltonian action of a Lie
group $G$ on $M$, with moment map $\Phi:\,M\ra\g^*$, such that this
action commutes with the action of $S^1$, there is a naturally induced
Hamiltonian $G$-action on $M_+$, which agrees with the given actions
and moment maps on $M_0$ and $\psi^{-1}(\R_{>0})$.
\end{proposition}
By reversing the $S^1$-action, one can also define a cut-space 
$M_-$, which is the union of $M_0$ and $\psi^{-1}(\R_{<0})$. 

The orbifold structure on $M_\pm$ depends only on the circle action on
$\psi^{-1}(0)$; replacing the given circle by a covering introduces
additional orbifold singularities.  This shows that cutting is a {\em
local} operation, that is, the construction only depends on the
existence of the circle action on $\psi^{-1}(0)$.  For example, one
can cut the two-torus $S^1\times S^1$, with its area form as a
symplectic form, along the circle $S^1\times\{1\}$ by using a locally
defined Hamiltonian for the $S^1\times\{1\}$-action. Note that in this
case the ``cut space'' $M_{cut}$ (which is the disjoint union $M_+\cup
M_-$ for global Hamiltonian $S^1$-actions) is connected, and is in
fact just the two-sphere.

Returning to the case of a globally defined Hamiltonian $S^1$-action,
suppose that we are given an $S^1\times G$-equivariant complex vector
bundle $E\ra M$. By pulling $E$ back to $M\times \C^-$, restricting to
$\ti{\psi}^{-1}(0)$ and taking the quotient
$$E_+:=(\text{pr}_1^*E|\ti{\psi}^{-1}(0))/S^1 $$
we obtain a $G$-vector bundle over $M_+$. There are canonical equivariant
isomorphisms 
$$ \iota_0^* E_+\cong E_0,\,\,\iota_{>0}^*\,E_+\cong E|\,\psi^{-1}(\R_{>0}). 
$$
In particular, one
obtains $G$-equivariant ``cut bundles'' $E_+\ra M_+$ and $E_-\ra M_-$
even if the $S^1$-action is only defined near $\psi^{-1}(0)$. (Of
course, the cut bundle will depend on the choice of a lift of the
$S^1$-action to $E$.)  To include the case that the cut space is
connected, we denote the cut bundle by $E_{cut}\ra M_{cut}$.

If $E=L$ is a pre-quantum line bundle, i.e. has a Hermitian 
metric and connection with curvature the symplectic form, then 
$L_+=(L\boxtimes\C)\qu S^1$ with induced 
metric and connection is a pre-quantum line bundle.  
The following result was proved in \cite{M2} and independently
(for the case of global $S^1$-actions) by Siye Wu.   

\begin{proposition}\label{CircleGluing}
Let $E\ra M$ be a $G$-equivariant vector bundle, such that the
$G$-action on $M$ preserves the symplectic form. Let $E_0\ra M_0$ be
the reduced bundle and $E_{cut}\ra M_{cut}$ the cut bundle, with
respect to a local Hamiltonian $S^1$-action.  The Riemann-Roch numbers
of the cut bundles satisfy the gluing rule
$$ RR(M,E)=RR(M_{cut},E_{cut})-RR(M_0,E_0).$$
\end{proposition}

By iterating the cutting operation, it is possible to obtain more
general cut-spaces. The idea is to take a ``suborbifold with 
corners'' $N\subset M$ with the property that the nullfoliation 
on the boundary faces of $N$ are orbifold-torus bundles; the 
cut space $M_N=N/\sim$ is  obtained by collapsing the boundary faces.   

We make this precise as follows.  Let $N \subset M$ be a closed subset
together with a finite collection of open subsets $U_j, \, j\in J$ of
$M$ equipped with Hamiltonian $S^1$-actions with moment maps
$\psi_j:\, U_j \ra \R$, and an open subset $V\subset M$ equipped with
a Hamiltonian action of a torus $H$, with moment map $\Psi:\,V\ra
\h^*$.
\begin{definition}
We will call a closed subset $N\subset M$ together with the collection
$(V,\Psi),\{(U_j,\psi_j)\}$ a {\bf sub-orbifold with reducible
corners} if the following holds:
\begin{enumerate}
\item Locally, near any point $n$, $N$ is given as the intersection 
$$\bigcap_{U_j\ni n}\psi^{-1}_j(\R_{\ge 0})\cap \Psi^{-1}(0).$$
\item On each non-empty intersection $U_{j_1}\cap\ldots\cap
      U_{j_k}\cap V$, the $S^1$-actions for the moment maps
      $\psi_{j_\nu}$ commute with each other and with $H$, and the
      corresponding action of $(S^1)^k\times H$ is locally free.
\end{enumerate}
\end{definition}
The {\bf boundary faces} of $N'\subset N$ are obtained by setting 
some of the ``boundary defining functions'' $\psi_j$ equal to zero. 
By including these into $\Psi$, one sees that every boundary face is
also a sub-orbifold with reducible corners.

By reducing with respect to the $H$-action and collapsing the boundary
facets $N\cap \psi_j^{-1}(0)$ as above, one obtains a cut space $M_N$
of dimension $\dim M-2\dim H$.  Given a complex vector bundle $E\ra V$
together with lifts of the local $H$- resp. $S^1$-actions so that the
lifts commute on the intersections $U_{j_1}\cap\ldots\cap U_{j_k}\cap
V$, one obtains cut bundles $E_N\ra M_N$.

\begin{remark}
\begin{enumerate}
\item 
Even though the nullfoliation of the boundary faces is intrinsically
defined, the cut space depends, as an orbifold, on the choice of the
local $S^1$-actions.
\item
In case $M$ is a Hamiltonian $G$-orbifold and $N$ a $G$-invariant
sub-orbifold with reducible corners, the cut space $M_N$ is a
Hamiltonian $G$-orbifold if the local $H$-and $S^1$-actions commute with
$G$.
\item If $E=L$ is a pre-quantum line bundle and if the local lifts of
the local circle actions satisfy the pre-quantum condition then $L_N$
is a pre-quantum line bundle.
\end{enumerate}
\end{remark}

The following Theorem shows how the Riemann-Roch number of the original 
bundle decomposes into Riemann-Roch numbers for cut bundles.
\begin{theorem} (Gluing Formula) \label{GluingThm}
Let $M$ be a compact Hamiltonian $G$-orbifold, $E\ra M$ a
$G$-equivariant complex orbifold vector bundle, and $\lie{N}=\{N\}$ a
collection of sub-orbifolds with reducible corners of $M$, together
with $G$-equivariant lifts to $E$ of the local $S^1$- and
$H$-actions. We assume that
\begin{enumerate}
\item The collection $\lie{N}$ covers $M$.
\item For each $N\in\lie{N}$, all boundary faces of $N$ are also in
$\lie{N}$.
\item The intersection of any two $N\in \lie{N}$ is either empty or is
a boundary face of each.
\end{enumerate}
Then the $G$-equivariant Riemann-Roch numbers of the cut spaces satisfy 
the gluing rule
\begin{equation} \label{GluingForm}
\RR(M,E)=\sum_{N\in \lie{N}}\,(-1)^{\text{codim} (N)} RR(M_N,E_N).
\end{equation}
\end{theorem}
For sub-orbifolds with reducible corners constructed from induced
toric maps for Hamiltonian actions of compact groups $G$ the above
gluing formula was proved in \cite{M2}; the general result will be
discussed in a separate paper \cite{MW}.

\subsection{Induced toric actions for compact groups}

Let $G$ be a compact connected Lie group, and $M$ a Hamiltonian
$G$-orbifold, with moment map $\Phi:\,M\ra \g^*$, and let
$\ti{\Phi}=q\circ \Phi$ be the corresponding induced toric map.
Recall from Section \ref{ToricSection} that for every open face
$\sig\subset \t^*_+$, there is a $G$-equivariant Hamiltonian action of
the center $Z(G_\sig)\subset G_\sig$ on $G\cdot
Y_\sig=\Phinv(U_\sig)$, with moment map the $\z(\g_\sig)^*$-component
of $\ti{\Phi}$.

A sub-orbifold with reducible corners $N \subset M$ can be constructed
as follows.  Let $Q\subset \t^*$ be a simple rational polytope, with
the property that for any two open faces $\sig$ of $\t^*_+$ and $F$ of
$Q$ with non-empty intersection, the tangent space to $F$ contains the
space perpendicular to $\sig$, i.e. the space $\g_\sig^0\cap \t^*$.
Examples of polytopes $Q$ with these properties are given in Figure
\ref{CutFig} below. We claim that the pre-image $N:=\ti{\Phi}^{-1}(Q)$
is a sub-orbifold with reducible corners for generic $Q$.

For each open face $F\subset Q$, let $T_F$ denote the the torus 
with Lie algebra the annihilator of the tangent space to $F$. 
In particular, let $H:=T_{\text{int}(Q)}$, and define
$$\Psi:=\text{pr}_{\h^*}(\ti{\Phi}-\mu)$$ 
where $\mu\in Q\cap\t^*_+$ and $\text{pr}_{\h^*}$ is the projection to 
$\h^*$. 

By assumption, $H\subset
Z(G_\sig)$ for every face $\sig$ that meets $Q$. It follows that
$\Psi$ is smooth on a $G$-invariant neighborhood $V$ of $N$ and
generates a $G$-equivariant Hamiltonian $H$-action. 

Similarly, given any codimension 1-face $F\subset Q$ with 
$F\cap\t^*_+\not=\emptyset$, let $v_F\in \Lambda$ be a lattice vector 
such that $\h\oplus \R\cdot v_F=\R\cdot(F-\mu_F)^0$ 
where $\mu_F\in F$. 
On a neighborhood $U_F$ of $\ti{\Phi}^{-1}(F)$, the map 
$$\psi_F=\l \ti{\Phi}-\mu_F,\,v_F\r: \, U_F\ra \R $$
is smooth and generates an $G\times H$-equivariant $S^1$-action.
Moreover, the $S^1$-actions on all overlaps $U_{F_1}\cap\ldots \cap
U_{F_k}$ commute.  The condition for the above data to define a
sub-orbifold with reducible corners is that for each open face
$F\subset Q$, the $T_F\subset T\subset G$-action on
$\Psi^{-1}(F\cap\t^*_+)$ is locally free. In this case, we call the
polytope $Q$ {\bf admissible} and denote the cut space by $M_Q=M_N$.
Note that if $Q=\{\mu\}$ with $\mu \in \text{int}(\t_+^*)$ then $Q$
is admissible if and only if $\mu$ is a regular value of $\Phi$, and
in this case $M_Q$ is just the reduced space $M_\mu$.

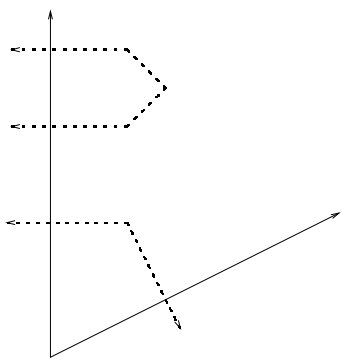
\begin{figure}[htb] 
\begin{center}

\setlength{\unitlength}{0.00033333in}
\begingroup\makeatletter\ifx\SetFigFont\undefined%
\gdef\SetFigFont#1#2#3#4#5{%
  \reset@font\fontsize{#1}{#2pt}%
  \fontfamily{#3}\fontseries{#4}\fontshape{#5}%
  \selectfont}%
\fi\endgroup%
{\renewcommand{\dashlinestretch}{30}
\begin{picture}(5199,5439)(0,-10)
\dashline{60.000}(1887,2112)(2712,462)
\path(2631.502,555.915)(2712.000,462.000)(2685.167,582.748)
\dashline{60.000}(1887,2112)(12,2112)
\path(132.000,2142.000)(12.000,2112.000)(132.000,2082.000)
\dashline{60.000}(2487,4212)(1887,4812)(87,4812)
\path(207.000,4842.000)(87.000,4812.000)(207.000,4782.000)
\dashline{60.000}(2487,4212)(1887,3612)(87,3612)
\path(207.000,3642.000)(87.000,3612.000)(207.000,3582.000)
\path(687,12)(687,5412)
\path(717.000,5292.000)(687.000,5412.000)(657.000,5292.000)
\path(687,12)(5187,2262)
\path(5093.085,2181.502)(5187.000,2262.000)(5066.252,2235.167)
\end{picture}
}

\caption{Examples of $Q$ for $G=SU(3)$ \label{CutFig}}
\end{center}
\end{figure}

If $E \ra M$ is a $G$-equivariant vector bundle then as explained in
Section \ref{ToricLift} there is a canonical lift of the 
$Z(G_\sig)\times G$-action on $\ti{\Phi}^{-1}(U_\sig)=G\cdot Y_\sig$. 
Therefore one also has a canonical lift of the local $H$- resp. $S^1$-actions, 
and a corresponding cut bundle $E_N$. 
However, if $E = L$ is a
$G$-equivariant pre-quantum bundle then this lift does  not in general
satisfy
the pre-quantum condition, and therefore 
the cut bundle $L_N$ is not pre-quantum.

However, if $Q$ is a lattice polytope then $\text{pr}_{\h^*}(\mu)$
lies in the lattice for $H\subset T$, and also $\l \mu_F,\nu_F\r\in
\Z$.  One therefore obtains local pre-quantum lifts by multiplying the
given $H$- resp. $S^1$-actions by the character $\exp(2\pi i\l
\text{pr}_{\h^*}(\mu),\xi\r)$ resp.  $\exp(2\pi i\l
\mu_F,\nu_F\r\phi)$. The cut bundle $L_Q$ for the modified actions is
a $G$-equivariant pre-quantum line bundle.  In case $Q$ has rational
vertices, choose any covering $\ti{T}\ra T$ such that $Q$ is a lattice
polytope in the refined lattice $\ti{\Lambda}^*\supset \Lambda^*$. The
local actions of the corresponding covers $\ti{H}\ra H$ and $S^1\ra
S^1$ admit pre-quantum lifts as before.  Notice again that passing to
a cover introduces additional orbifold singularities.

A family $\lie{N}$ satisfying the hypotheses of Theorem
\ref{GluingThm} can be constructed as follows.  Let $M$ be a
compact Hamiltonian $G$-orbifold, and $\mu \in \text{int}(\t_+^*)$.
Recall that the {\bf dual cone} $C_\sig$ to $\t^*_+$ at an open face
$\sig\subset \t^*_+$ is defined to be the set of all vectors $v \in
\t^*$ such that the minimum of the inner product $(v,x)$ for $x \in
\t^*_+$ is achieved by $x \in \sig$.  For each face $\tau \subset
\ol{\sig}$ let $Q_{\sig,\tau}$ denote the polytope $\mu+\ol{\tau} -
C_\sig$, and ${\mathcal{Q}}$ the collection of these polytopes. For
generic choices of $\mu$, all of these polytopes are admissible.

\begin{figure}[htb]
\begin{center}
\setlength{\unitlength}{0.00033333in}
\begingroup\makeatletter\ifx\SetFigFont\undefined%
\gdef\SetFigFont#1#2#3#4#5{%
  \reset@font\fontsize{#1}{#2pt}%
  \fontfamily{#3}\fontseries{#4}\fontshape{#5}%
  \selectfont}%
\fi\endgroup%
{\renewcommand{\dashlinestretch}{30}
\begin{picture}(5124,5739)(0,-10)
\path(612,312)(5112,2562)
\path(5018.085,2481.502)(5112.000,2562.000)(4991.252,2535.167)
\path(612,312)(612,5712)
\path(642.000,5592.000)(612.000,5712.000)(582.000,5592.000)
\dashline{60.000}(1212,5712)(1212,1212)
\path(1242.000,5592.000)(1212.000,5712.000)(1182.000,5592.000)
\dashline{60.000}(1287,1212)(5112,3312)
\path(5021.248,3227.952)(5112.000,3312.000)(4992.373,3280.546)
\dashline{60.000}(1212,1212)(1812,12)
\path(1731.502,105.915)(1812.000,12.000)(1785.167,132.748)
\dashline{60.000}(1212,1212)(12,1212)
\path(132.000,1242.000)(12.000,1212.000)(132.000,1182.000)
\end{picture}
}
\caption{The collection $\mathcal{Q}$ for $G=SU(3)$}
\end{center}
\end{figure}
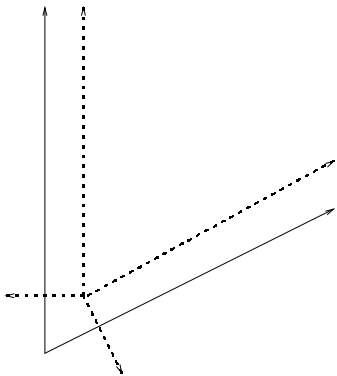

Next consider the case that $M$ is a Hamiltonian $G \times
G$-orbifold, with moment map $(\Phi_+,\Phi_-)$ and assume that $0$ is
a regular value for the diagonal $G$-action.  Let $M_0 = M \qu G$ be
the reduced space and $\ti{\Phi}_0:\,M_0\ra \t^*_+$ the residual toric
moment map induced from the map $\ti{\Phi}_+=q\circ \Phi_+$ on the
first factor.  Given any polytope $Q$ as above the subset $N :=
\ti{\Phi}^{-1}_0(Q)$ is a sub-orbifold with reducible corners.

If $E \ra M$ is a $G\times G$-equivariant vector bundle, the quotient
bundle $E_0$ has a canonical lift of the local $Z(G_\sig)$-actions on
$M_0$, so that there is a canonical cut bundle.  The same argument as
above shows that if $E=L$ is a pre-quantum line bundle and $Q$ is a
rational polytope, one can modify the given lift of the local actions
in such a way that the cut bundle is a pre-quantum line bundle.

As an application, we prove:

\begin{proposition}\label{BabyVerlinde}  
Let $G$ be a compact connected Lie group and $M$ a Hamiltonian
$G\times G$-orbifold with proper moment map $\Phi=(\Phi_+,\Phi_-)$.
Let $L \ra M$ be a $G\times G$-equivariant pre-quantum line bundle.
If the symplectic quotient $M_0=M\qu G$ by the diagonal action is
compact, the Riemann-Roch number of the reduced bundle $L_0=L\qu G$ is
given by
\begin{equation}\label{DiagonalReduction}
\RR(M_0,L_0)=
\sum_{\mu \in \Lambda^*_+} \,\RR(M_{\mu,*\mu},L_{\mu,*\mu}) . 
\end{equation}
(Here all $\RR$-numbers are defined by partial desingularization 
if necessary.) 
\end{proposition}

\begin{proof}
If $M$ is compact, both sides are equal, by the ``quantization
commutes with reduction'' Theorem \ref{QuantizationReduction}.  If $M$
is not compact, let $\ti{\Phi}_0:\,M_0 \ra \t^*$ denote the 
residual toric moment map
induced by $\ti{\Phi}_+$.  Its image is given by
$$\ti{\Phi}_0(M_0) = \ti{\Phi}_+(M) \cap * \ti{\Phi}_-(M).$$
Choose a compact admissible rational polytope $Q\subset \t^* \times
\t^*$, such that $ \ti{\Phi}_0(M_0) \times *\ti{\Phi}_0(M_0) $ is
contained in the relative interior of $Q\cap ( \t^*_{+} \times
\t^*_{+})$.  We can assume that $Q$ has rational vertices, so that the
cut bundle $L_Q \ra M_Q$ is a $G$-equivariant pre-quantum bundle.
Since $M_Q$ is compact and $L\qu G=L_Q\qu G$, this reduces the proof
to the case where $M$ is compact.
\end{proof}

\subsection{Induced toric actions for loop groups}

Now let $G$ be a compact, connected, simply-connected Lie group, and
$\Alc$ the corresponding fundamental alcove.  Let $M$ be a Hamiltonian
$L(G\times G)$-manifold with proper moment map
$\Phi=(\Phi_{+},\Phi_{-})$ (at level $\lev=+1$).  We assume that the
anti-diagonal $LG$-action is locally free on its zero level set, so
that the symplectic quotient $M_0=M\qu LG$ is a finite dimensional,
compact symplectic orbifold.  Let $\ti{\Phi}_0:\,M_0\ra \Alc$ denote
the residual toric moment map induced from $q\circ \Phi_+$.  For each
face $\sig \subset \Alc$, the open set $\ti{\Phi}_0^{-1}(U_\sig)
\subset M_0$ carries a Hamiltonian action of $Z((LG)_\sig)$, with moment
map given by $\ti{\Phi}_0$ followed by projection to $\z((L \g)_\sig)^*$.

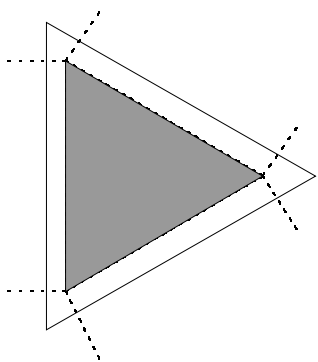
\begin{figure}[htb] 
\begin{center}
\setlength{\unitlength}{0.00033333in}
\begingroup\makeatletter\ifx\SetFigFont\undefined%
\gdef\SetFigFont#1#2#3#4#5{%
  \reset@font\fontsize{#1}{#2pt}%
  \fontfamily{#3}\fontseries{#4}\fontshape{#5}%
  \selectfont}%
\fi\endgroup%
{\renewcommand{\dashlinestretch}{30}
\begin{picture}(4824,5439)(0,-10)
\path(612,5262)(612,462)(4812,2862)
	(612,5262)(612,5262)
\dashline{60.000}(912,1062)(1437,12)
\dashline{60.000}(912,1062)(12,1062)
\dashline{60.000}(912,4662)(12,4662)
\dashline{60.000}(912,4662)(1437,5412)
\dashline{60.000}(3987,2862)(4512,3612)
\dashline{60.000}(3987,2862)(4512,2037)
\texture{55888888 88555555 5522a222 a2555555 55888888 88555555 552a2a2a 2a555555 
	55888888 88555555 55a222a2 22555555 55888888 88555555 552a2a2a 2a555555 
	55888888 88555555 5522a222 a2555555 55888888 88555555 552a2a2a 2a555555 
	55888888 88555555 55a222a2 22555555 55888888 88555555 552a2a2a 2a555555 }
\dashline{60.000}(912,4662)(912,1062)(3987,2862)
	(912,4662)(912,4662)
\shade\path(912,4662)(912,1062)(3987,2862)
	(912,4662)(912,4662)
\end{picture}
}
\caption{Cutting the fundamental alcove for $G=SU(3)$ \label{AlcoveFig}}

\end{center}
\end{figure}

Let
$$\Alc_\eps:=(1-\eps) \Alc + \eps \mu $$
for some small $\eps \in \R_{>0}$ and $\mu \in \text{int}(\Alc)$.
This is the shaded region in Figure \ref{AlcoveFig}.

For each face $\sig$ of $\Alc$, let $\sig_\eps = (1-\eps) \sig + \eps
\mu$ denote the corresponding face of $\Alc_\eps$.  For any face $\tau
\subset \ol{\sig}$ let $Q_{\tau,\sig} = \ol{\tau_\eps} - C_\sig$.  Let
${\mathcal{Q}}$ be the collection of polytopes $Q_{\tau,\sig}$ (see
Figure \ref{AlcoveFig}).  For generic choices of $\eps,\mu$, each $Q
\in \mathcal{Q}$ gives rise to a sub-orbifold with reducible corners
$N=\ti{\Phi}_0^{-1}(Q)$ and cut space $(M_0)_Q:=(M_0)_N$.

Also, for every $\widehat{LG^2}$-equivariant vector 
bundle $E\ra M$, any
lift of the local $Z((LG)_\sig)$-actions to $E_0=E\qu K$ defines a
cut-bundle $(E_0)_Q$. Theorem \ref{GluingThm} gives an expression for
$\RR(M_0,E_0)$ in terms of the Riemann-Roch numbers
$\RR((M_0)_Q,(E_0)_Q)$.

If $L\ra M$ is a $\widehat{LG^2}$-equivariant pre-quantum line bundle,
and $\eps \in \Q, \mu \in \Lambda^* \otimes_\Z \Q$ are chosen so that
each $Q$ is a polytope with rational vertices, one can pass to a cover
$\ti{T}\ra T$ and arrange the lifts in such a way that one obtains a
pre-quantum cut-bundle $(L_0)_Q$.

\section{Proof of Theorem \ref{MainResult}}

As in the statement of the theorem, let $G$ be compact, connected and
simply connected, and $\Alc$ its fundamental alcove.  Let $M$ be a
Hamiltonian $L(G\times G)$-manifold with proper moment map at level
$\levi \in \N$ and pre-quantum line bundle $L \ra M$, and $L\qu LG \ra
M\qu LG$ the symplectic reduction with respect to the anti-diagonal
$LG$-action.

The basic idea of the proof of \ref{MainResult} is to cut the
symplectic quotient $M\qu LG$ into pieces, all of which are obtained
by diagonal reduction of a finite-dimensional symplectic manifold by a
Hamiltonian action of a compact group.  The result then follows by
applying the ``quantization commutes with reduction'' Theorem for
compact groups, together with the gluing formula.  If the quotient
$M\qu LG$ is singular, we have to combine this idea with partial
desingularization.

Let us first assume that $0$ is a regular value for the anti-diagonal
action, so that $M\qu LG$ is a compact symplectic orbifold, and $L\qu
LG \ra M\qu LG$ a pre-quantum line bundle.  Let $\mathcal{Q}$ be a
collection of admissible rational polytopes in $\Alc$ defined in the
previous subsection.  By the Gluing Formula, Equation
(\ref{GluingForm}),
$$ \RR(M\qu LG,L\qu LG) = \sum_{Q \in\mathcal{Q} } \,(-1)^{\text{codim}(Q)}
\RR((M\qu LG)_Q,(L\qu LG)_Q). $$
Each cut space $(M\qu LG)_Q$ can be written as a reduction in finite
dimensions:  Let $\sig \subset \Alc$ be an open face
such that $Q \cap \Alc$ is contained in $U_\sig$, and let
$Y_{\sig,-\sig}$ denote the corresponding cross-section of $M$.
The restriction to $Y_{\sig,-\sig}$ of  $\ti{\Phi}_+$
defines an induced toric action on $Y_{\sig,-\sig}$; let $(Y_{\sig,-\sig})_Q$ 
be the corresponding cut space.  Then
$$ (M\qu LG)_Q = (Y_{\sig,-\sig})_Q \qu {K}_\sig , \ \ \ \ (L\qu LG)_Q =
  (L | Y_{\sig,-\sig})_Q \qu {K}_\sig.$$
By Proposition \ref{BabyVerlinde} applied to $(L | Y_{\sig,-\sig})_Q$
we obtain that
$$ \RR((M\qu LG)_Q, (L\qu LG)_Q) = \sum_{\mu \in Q \cap \Lambda^*_\levi}
\RR(M_{\mu,*\mu},L_{\mu,*\mu})$$
where $\Lambda^*_\levi = \Lambda^* \cap \levi \Alc$ is the set of
dominant weights for $LG$ at level $\levi$.  Finally, applying the
Euler identity $\sum_{Q \ni \mu} (-1)^{\text{codim}(Q)} = 1$,
\begin{eqnarray*}
\RR(M\qu LG,L\qu LG) &=& \sum_{Q\in\mathcal{Q}}(-1)^{\text{codim}(Q)}
\sum_{\mu \in Q \cap \Lambda^*_\levi}
\RR(M_{\mu,*\mu},L_{\mu,*\mu})\\
&=& \sum_{\mu \in \Lambda^*_\levi}
\RR(M_{\mu,*\mu},L_{\mu,*\mu})
\end{eqnarray*}
which proves the result in this case.

In the case where $0$ is not a regular value, we use a
desingularization by ``adding a puncture'', as explained in section
\ref{DesingSec}. Let 
$$ M^{({2})}:=M\times
\M(\Sig^3_0)\qu LG$$
and
$$M^{(1)}= M^{({2})}\qu LG=M\times
\M(\Sig^3_0)\qu LG\times LG,$$
equipped with their pre-quantum bundles. Correspondingly, we will use
cross-sections $Y^{(1)}_{\{0\}}$ for $M^{(1)}$,
$Y_{\{0\},\sig,-\sig}^{({2})}$ for $M^{({2})}$, and $Y_{\sig,-\sig}$ for
$M$. To compactify $Y^{(1)}_{\{0\}}$, we choose an admissible polytope
$Q_1$ with rational vertices with $0\in\text{int}(Q_1)$, and use the
cut space $(Y^{(1)}_{\{0\}})_{Q_1}$.  We also choose an admissible
family of polytopes $\mathcal{Q}_2$ to cut $Y_{\sig,-\sig}$ and
$(Y_{\{0\},\sig,-\sig}^{({2})})_{Q_1}$ into pieces. Then
$$ (Y_{\{0\},\sig,-\sig}^{({2})})_{Q_1\times Q_2}\qu G\cong
(Y_{\sig,-\sig})_{Q_2},\,\,\,\,\,\,\, (Y_{\{0\},\sig,-\sig}^{({2})})_{Q_1\times
Q_2}\qu (LG)_\sig=(Y_{\{0\}}^{(1)})_{Q_1\times Q_2}.
$$
Using the gluing formula, together with repeated application of
``quantization commutes with reduction'' we compute (omitting the line
bundles from the notation):
$$
\begin{array}{llll} 
\RR(M\qu LG)&=& \RR((Y_{\{0\}}^{(1)})_{Q_1}\qu G)
=\RR((Y_{\{0\}}^{(1)})_{Q_1})^G & \text{(Theorem
\ref{QuantizationReduction})} \\ &&& \\ 
&=&\sum_{Q_2\in\mathcal{Q}_2}
(-1)^{\text{codim}(Q_2)}\,\RR((Y_{\{0\}}^{(1)})_{Q_1\times Q_2})^{G} &
\text{(Theorem \ref{GluingThm})} \\ &&& \\
&=&\sum_{Q_2\in\mathcal{Q}_2} (-1)^{\text{codim}(Q_2)}\,
\RR((Y_{\{0\},\sig,-\sig}^{({2})})_{Q_1\times Q_2})^{(LG)_\sig\times G} &
\text{(Theorem \ref{QuantizationReduction})} \\ &&& \\
&=&\sum_{Q_2\in\mathcal{Q}_2}
(-1)^{\text{codim}(Q_2)}\,\RR((Y_{\sig,-\sig})_{Q_2})^{(LG)_\sig} &
\text{(Theorem \ref{QuantizationReduction})} \\ &&& \\ &=&
\sum_{Q_2\in\mathcal{Q}_2} (-1)^{\text{codim}(Q_2)}\sum_{\mu\in
Q_2\cap \Lambda^*_\levi} \RR(M_{\mu,*\mu})& \text{(Theorem
\ref{BabyVerlinde})} \\ &&& \\ &=& \sum_{\mu\in \Lambda^*_\levi}\RR(M_{\mu,*\mu}) & \text{(Euler Identity)}.
\end{array} 
$$

\begin{appendix}

\section{The gauge-theoretic construction of $\mathcal{M}(\Sigma)$}
\label{YangMills}

Let $\Sig$ be a compact Riemann surface of genus $g$ with $b$ boundary
components, and $G$ a connected, simply connected compact Lie group.
Let $\A(\Sig)\cong\Omega^1(\Sig,\g)$ be the space of $G$-connections
of Sobolev class $s>\f{1}{2}$.  Recall that for any manifold $X$ with 
boundary, the restriction map $C^\infty(X)\to C^\infty(\p X)$ to the boundary 
extends for $t >\f{1}{2}$ to a continuous surjection   
$$H_{(t)}(X)\ra H_{(t-\f{1}{2})}(\p X)$$ 
(see e.g. \cite{BW}, chapter
11). Hence there is a surjective map from $\A(\Sig)$ to the space 
$\A(\p\Sig)\cong\Omega^1(\p\Sig,\g)$ of connections 
over $\p\Sig$ of Sobolev class $s-\f{1}{2}$. 
The space $\G(\Sig)$ of maps from $\Sig$ to $G$ of Sobolev class $s+1$ consists
of continuous maps, and is a Banach Lie group  acting
smoothly on the space $\A(\Sig)$. 
Restriction to the
boundary gives a surjective map from $\G(\Sig)$  to the space 
$\G(\partial \Sig)$ of 
maps $\p\Sig\to G$ of Sobolev class $s+\f{1}{2}$, and
 as before we  define $\G_\p(\Sig)$ to be the kernel.

Let us fix a Riemannian metric on $\Sig$, and parametrizations 
$B_j\cong S^1$ of the boundary components compatible with the metric 
and orientation. As before we identify $\G(B_i)\cong LG$ and 
$\A(B_i)\cong L\g^*$. If we take $\Omega^2(\Sig,\g)$ 
to consist of two-forms of Sobolev-class $s-1$, the moment map for 
the $\G(\Sig)$-action becomes a continuous map
$$\A(\Sig)\ra \Omega^2(\Sig,\g)\oplus \Omega^1(\p\Sig,\g),\,A\mapsto (F_A,\,
\iota^*A).$$

An atlas for the moduli space 
$\M(\Sig) = \A_F(\Sig)/\G_\partial(\Sig)$ for $b>0$ is
constructed, as in the case without boundary, from local slices 
for the gauge group action \cite{D,DK,AHS}.  
Let $A \in \A(\Sig)$ be any connection.  By the implicit function
theorem, any connection $A + a$, with $a$ small, can be gauge
transformed by a unique element $g\in \G_\partial(\Sig)$ into Coulomb
gauge with respect to $A$, that is, so that
\begin{equation}
d_A^*(g\cdot(A+a)-A)=0.   \label{CoulombGauge}
\end{equation}
In other words, a neighborhood of $A$ in $A+\text{ker}(\d_{A}^*)$ is a
slice for the $\G_\p(\Sig)$-action on $\A(\Sig)$.  For any $A \in
\A_F(\Sig)$ one defines a local moduli space to be a neighborhood of
$A$ in $\A_F(\Sig)$ of the form
$$ V_A \subset \{
A+a \in \Omega^1(\Sig,\g)|\, \ F_{A+a} = 0,\, d_A^*a=0\}. $$ 
Using the implicit function theorem again,
one shows that if $V_A$ is taken sufficiently small, 
$V_A \subset \A_F(\Sig)$ is a smooth Banach submanifold
locally homeomorphic to its tangent space at $0$
$$ T_0(V_A) = \{
a \in \Omega^1(\Sig,\g)|\,\d_A\,a=0,\,\d_A^*\,a=0\}. $$ 
The sets $V_A$
together with the coordinate mappings $V_A \ra T_0(V_A)$ give an atlas
for $\M(\Sig)$.  Since the tangent space $T_0(V_A)$ is invariant under
the Hodge $*$-operator, it is a complex (and therefore symplectic)
subspace of $T_A(\A(\Sig))$ and this shows that $\M(\Sig)$ is a Banach
K\"ahler manifold.

Similarly, the products $V_A\times \C\cong T_0(V_A)\times \C$ are 
slices for the $\G_\p(\Sig)$-action on $\A_F(\Sig)\times\C$, and 
give local bundle charts for the line bundle $L(\Sig)$.

\begin{remark}
For a Riemann surface $\Sig$ without boundary, the moduli space
$\M(\Sig)$ is singular in general. Its smooth part $\M(\Sig)_{smooth}$
is given by gauge equivalence classes of irreducible connections, i.e.
those with trivial stabilizer. The above proof applies to the smooth
part and shows that it is a symplectic manifold.
\end{remark}

We are now in position to prove Theorem \ref{GluingEqualsReduction}. 
\label{theproof}

\begin{theorem}
\label{GluingEqualsReduction2}
Let ${\Sig}$ be a compact oriented Riemann surface with $b>0$ boundary
components.  Suppose $\Sig$ is obtained from a compact oriented
Riemann surface $\hat{\Sig}$ by gluing along two boundary components
$B_\pm \subset \partial \hat{\Sig}$.  Let $\pi:\hat{\Sig}\to \Sig$ be
the gluing map. Then the map
\begin{equation} \label{Homeo}
\M(\Sig) \ra \M(\hat{\Sig}) \qu LG,\,\,
[A] \mapsto LG\cdot [\pi^*A]\end{equation}
is an $LG^b$-equivariant symplectic diffeomorphism.  Similarly, there
is an $\widehat{LG^b}$-equivariant isomorphism of line bundles
$$ L(\Sig) = L(\hat{\Sig}) \qu LG.$$
\end{theorem}

\begin{proof}
Choose a parametrization $B_\pm \cong S^1$ as in Section
\ref{SubsectionGluingEqualsReduction}, and consider the corresponding
anti-diagonal $LG$-action on $\M(\hat{\Sig})$.  Let $A\in\A_F(\Sig)$ be
any fixed connection, and $V_A\subset\A_F(\Sig)$ the corresponding
local moduli space.  Define a map
$$ \psi:\,\,V_A\times LG \ra \M(\hat{\Sig}),\,(A+a,h)\mapsto
 {h}\cdot[\pi^*(A +a)].$$
Notice that this map is $LG$-equivariant and that its image is
contained in the zero level set of $\M(\hat{\Sig})$. To show that
$\psi$ is a diffeomorphism from a small neighborhood of $(0,e)$ onto
its image, we can assume after acting by a suitable element of
$\G(\Sig)$ that $A$ is smooth.  Since the map in (\ref{Homeo}) is a
homeomorphism it suffices to prove that $\psi$ is an immersion, that is,
the tangent map to $\psi$ at $(0,e)$ is injective and has closed
image.  In terms of the local moduli space $V_{\pi^*A}$,
the map $\psi$ is given as
\begin{equation} \label{Label}
 \psi:\,\,V_A\times LG \ra V_{\pi^*A},\,(A+a,h)\mapsto
\ti{h}\cdot\pi^*(A +a)\end{equation}
where $\ti{h}\in\G(\hat{\Sig})$ is the unique gauge transformation
such that $\ti{h}|\,B_\pm=h$ and $\ti{h}|\p\Sig-(B_+\cup B_-)=1$, 
and such the right-hand side of (\ref{Label})
lies in $V_{\pi^* A}$.  The tangent map is given by
$$ \d_{(0,e)}\psi(b,\eta)= \pi^*b - \d_{\pi^*A} \ti{\eta} $$
where $\ti{\eta} \in\Omega^0 (\hat{\Sig},\g)$ is the unique solution of the
Dirichlet problem
$$\d_{\pi^* A}^*\d_{\pi^* A} \ti{\eta}=0,\  
\iota^*_\pm \ti{\eta} = \eta |\, C,\  
\ti{\eta} |\,(\p\hat{\Sig}-B_+-B_-)=0.$$
Suppose that $(a,\eta) \in \ker d_{(0,e)} \psi$. Since $\ti{\eta}$ is
continuous, and since the restrictions of $\ti{\eta}$ to $B_\pm$ are
equal, $\ti{\eta} = \pi^* \zeta$ is the pullback of a continuous
function ${\zeta} \in C^0(\Sig,\g)$. Then $d_A \zeta = a$ in the sense
of distributions, which by ellipticity of $\d_A$ on $0$-forms implies
that the $\zeta$ is smooth.  Since $d_A^* d_A \zeta = 0$ and $\zeta
\vert {\p \Sig} = 0$ it follows that ${\zeta} = 0$ and hence $a = 0$,
as required.  This shows the first claim.  The image of $\d \psi$ is
closed by ellipticity of $\d_{\pi^*A}$.
\end{proof}

\begin{remark}
In the case that $\Sig$ has empty boundary so that $\M(\Sig)$ is
possibly singular, the above proof gives a symplectomorphism between
the smooth parts of $\M(\hat{\Sig})\qu LG$ and $\M(\Sig)$. In order 
to extend this to an isomorphism of stratified symplectic spaces 
in the sense of Sjamaar-Lerman \cite{SL}, it would be necessary to define 
a Poisson algebra of ``smooth'' functions on $\M(\Sig)$. In this paper, 
we do this by representing $\M(\Sig)$ as a symplectic quotient  
$\M(\Sig)=\M(\Sig')\qu LG$, where $\Sig'$ denotes $\Sig$ 
with one puncture 
added, and use symplectic cross-sections to replace $\M(\Sig')\qu LG$ 
by a finite dimensional symplectic quotient. Similarly, the singular 
symplectic quotient $M(\hat{\Sig})\qu LG$ is defined by adding a 
puncture. 
\end{remark}

\section{The action of $\widehat{\mathcal{G}(\Sig)}$ on 
$\A(\Sig)\times \C $} Let $\Sig$ be a compact, connected oriented
Riemann surface with boundary, $G$ a connected, simply connected
compact Lie group and $\A(\Sig)$ the space of flat $G$-connections (of
Sobolev class $s>\f{1}{2}$).  In this appendix, we will explain,
following Mickelsson \cite{Mi}, Ramadas-Singer-Weitsman \cite{RSW},
and Witten \cite{W2} how to obtain an explicit construction of the
canonical central extension $\widehat{\G(\Sig)}$ of the gauge group
and of its action on the pre-quantum line bundle $\A(\Sig)\times \C$.

Let us briefly recall some facts about pre-quantization. Let
$(M,\omega)$ be a connected symplectic manifold, and $L\ra M$ a
pre-quantum line bundle with pre-quantum connection $\nabla$.  By
Kostant's theorem \cite{Ko}, the group $\text{Aut}_\theta(L)$ of
connection preserving Hermitian bundle automorphisms is a central
extension by $S^1$ of the group $\text{Diff}_\omega(M)$ of
symplectomorphisms of $M$. Hence, given a Lie group $K$ and a
symplectic action $K\ra \text{Diff}_\omega(M)$ one obtains a central
extension
\begin{equation}\label{CentralExtension}
1 \ra S^1\ra \widehat{K} \ra K \ra 1, 
\end{equation}
together with an action $\widehat{K}\ra  \text{Aut}_\theta(L)$. 
If the action of $K$ is in fact Hamiltonian, with equivariant 
moment map 
$\Phi:\,M\ra \k^*$, the formula \ref{LiftFormula}
defines a Lie algebra splitting $\hat{\k}\cong\R\oplus \k$
or equivalently an invariant flat connection on the bundle 
(\ref{CentralExtension}). 
If $\pi_0(K)=\{1\}$ and if the holonomy of this connection  
is trivial, one obtains in this way a lift $K\ra \text{Aut}_\theta(L)$. 

Let us apply these results to the case $M=\A(\Sig), K=\G(\Sig)$ and
$L= \A(\Sig)\times \C$, equipped with the pre-quantum connection
$\theta_A(a)=\f{1}{2} \int_\Sig\, a\stackrel{ _\cdot}{\wedge}
A $. Using the fact that $\pi_0(G),\,\pi_1(G),\,\pi_2(G)$ are all
trivial, one can show that the gauge group $\G(\Sig)$ (and also
$\G_\p(\Sig)$) is connected:
$$ \pi_0(\G(\Sig))=\pi_0(\G_\p(\Sig)) =\{1\}. $$
However, the fundamental groups of $\G(\Sig),\,\G_\p(\Sig)$ are
non-trivial in general.

Recall now the description of the central extension $\widehat{\G(\Sig)}$ 
given in section \ref{ConstructionModuliSpace}. On the Lie algebra 
level, $\text{Lie}(\widehat{\G(\Sig)})$ is the central extension 
by $\R$ of $\text{Lie}(\G(\Sig))$ defined by the cocycle 
$$ (\xi_1,\xi_2)\mapsto \f{1}{2\pi} \int_\Sig \d\xi_2 \stackrel{ _\cdot}{\wedge}
\d\xi_2=  \f{1}{2\pi}\int_{\p\Sig} \xi_1\d\xi_2,$$
that is, $\text{Lie}(\widehat{\G(\Sig)})=\text{Lie}(\G(\Sig))\oplus\R$
with bracket 
$$ [(\xi_1,t_1),(\xi_2,t_2)]=\Big( [\xi_1,\xi_2],\, \f{1}{2\pi}\int_{\p\Sig}
\xi_1\d\xi_2\Big). $$
One checks that the pre-quantum lift (\ref{LiftFormula}) corresponds
to the inclusion of $\text{Lie}(\G(\Sig))$ in
$\text{Lie}(\widehat{\G(\Sig)})$ as the first summand.  Notice however
that since the moment map for the $\G(\Sig)$-action is equivariant
only over the subgroup $\G_\p(\Sig)$, the connection on
$\widehat{\G(\Sig)}$ obtained in this way is flat only over
$\G_\p(\Sig)$.  We will now check that this flat connection has
trivial holonomy, or equivalently that the restriction to
$\G_\p(\Sig)$ of the cocycle $c$ in Equation (\ref{cocycle}) is a
coboundary.

Suppose first that $\p\Sig=\emptyset$, so that $\G_\p(\Sig)=\G(\Sig)$. 
Let 
$\lambda\in\Omega^3(G)$ be the closed bi-invariant 3-form on $G$ 
given by 
$$\lambda=\f{1}{6} g^{-1}\,\d g \stackrel{ _\cdot}{\wedge} [g^{-1}\,\d
g,\,g^{-1}\,\d g].$$
Then $\f{1}{8\pi^2}\lambda$ (using the normalized inner product
defined in Section \ref{SubsectionGluingEqualsReduction}) defines an
integral cohomology class.  Let $B$ be any compact oriented 3-manifold
with boundary $\partial B={\Sigma}$ and $\hat{g}\in \text{Map}(B,G)$
an extension of $g$ of Sobolev class $s+\f{3}{2}$.  The integral over
$B$ of $\f{1}{8\pi^2}\hat{g}^*(\lambda)$ is independent of $B$ mod
$\Z$, and therefore
$$\Gamma:\,\G(\Sigma)\ra S^1,\,\,g\mapsto \exp\big(\f{i}{4\pi}\int_B 
\hat{g}^*(\lambda)\big)$$ 
a well-defined smooth map. One verifies the coboundary property 
\begin{equation}\label{CocycleProperty}
\Gamma(g_1\,g_2)=\Gamma(g_1)\,\Gamma(g_2)\,c(g_1,g_2).
\end{equation}
The group homomorphism $\G(\Sig)\ra \widehat{\G(\Sig)},\,g\mapsto
(g,\Gamma(g))$ defines the trivialization.
 
In the case $\p\Sig\not=\emptyset$, consider the Riemann surface
$\ol{\Sigma}$ without boundary obtained from $\Sigma$ by ``capping
off'' the boundary components, and let $\ol{\Gamma}:\,\G(\ol{\Sig})\ra
S^1$ be defined as before.  Let $\G_c(\Sig) \subset \G({\Sig})$ be the
subgroup of gauge transformations in the kernel of the restriction map
$\G(\ol{\Sig})\ra \G(\ol{\Sig}-\Sig)$. Then the restriction $\Gamma$
of $\ol{\Gamma}$ to $\G_c(\Sig)$ satisfies (\ref{CocycleProperty}) and
therefore shows that the extension is trivial over $\G_c(\Sig)$.  Now
consider $\G_\p(\Sig)\supset\G_c(\Sig)$. One sees easily that any loop
in $\G_\p(\Sig)$ can be deformed to a loop in $\G_c(\Sig)$, i.e. that
the natural map $\pi_1(\G_c(\Sig))\ra \pi_1(\G_\p(\Sig))$ is
surjective. It follows that the central extension is trivial over
$\G_\p(\Sig)$ as well.

The quotient $\widehat{\G(\p\Sig)}=\widehat{\G(\Sig)}/\G_\p(\Sig)$ is
the (unique) central extension of $\G(\p\Sig)$ defined by the Lie
algebra cocycle $(\xi_1,\xi_2)\mapsto \f{1}{2\pi}\int_{\p\Sig}\xi_1\cdot
\d\xi_2$. It follows that $\widehat{\G(\p\Sig)}$ is just the basic
central extension $\widehat{LG^b}\to LG^b$.  (In fact, the above
construction with $\Sig$ the two-disk is precisely Mickelsson's
construction \cite{Mi} of $\widehat{LG}$.)

\end{appendix}

\end{document}